\documentclass[a4paper, 10pt]{article}

\usepackage[a4paper,left=3cm,right=3cm,top=2.5cm,bottom=2.5cm]{geometry}
\usepackage{hyperref}
\usepackage{subfig}
\usepackage{pgfplots}
\usepackage{pgfplotstable}
\usepackage{mathtools}
\usepackage{multicol}
\usepackage{comment}
\usepackage{booktabs}
\pgfplotsset{compat=1.5}
\usepackage{amssymb}
\usepackage{url}
\usepackage{bm}

\usepackage{pdfsync}
\usepackage{float}
\usepackage{tabularx}
\usepackage{enumerate}
\usepackage{array}
\usepackage{xspace}
\usepackage{tikz}
\usepackage{tikz-cd}
\usepackage{tikzsymbols}
\usetikzlibrary{calc, trees, positioning,arrows, chains, shapes.geometric,
    decorations.pathreplacing, decorations.pathmorphing, shapes,
    matrix, shapes.symbols, decorations.markings, patterns,fit}

\usepackage{siunitx}
\usepackage{authblk}

\usepackage{mathrsfs}
\usepackage{color, colortbl}
\usepackage{multirow}

\usepackage{amsthm}
\usepackage{amsmath}
\usepackage{stmaryrd}
\usepackage{algorithm}
\usepackage{algpseudocode}
\usepackage{graphicx}

\usepackage[noabbrev]{cleveref}
\usepackage[normalem]{ulem}

\newcommand{\mupar}{\ensuremath{\boldsymbol{\mu}}}
\newcommand{\xpar}{\ensuremath{\mathbf{x}}}
\newcommand{\Xpar}{\ensuremath{\mathbf{X}}}
\newcommand{\Dpar}{\ensuremath{\mathcal{D}}}
\newcommand{\ypar}{\ensuremath{\mathbf{y}}}
\newcommand{\Ypar}{\ensuremath{\mathbf{Y}}}

\DeclareMathOperator*{\argmin}{\arg\!\min}
\DeclareMathOperator*{\argmax}{\arg\!\max}
\newcolumntype{C}[1]{>{\centering\arraybackslash}m{#1}}

\definecolor{Gray}{gray}{0.9}

\begin{document}

\title{Data-driven parameterization refinement for the structural optimization of cruise ship hulls}

\author[a]{Lorenzo~Fabris\footnote{lorenzo.fabris@sissa.it}}
\author[b]{Marco~Tezzele\footnote{marco.tezzele@emory.edu}}
\author[c]{Ciro~Busiello\footnote{ciro.busiello@fincantieri.it}}
\author[c]{Mauro~Sicchiero\footnote{mauro.sicchiero@fincantieri.it}}
\author[a]{Gianluigi~Rozza\footnote{gianluigi.rozza@sissa.it}}

\affil[a]{Mathematics Area, mathLab, SISSA, Scuola Internazionale Superiore di Studi Avanzati, Trieste, Italy}
\affil[b]{Department of Mathematics, Emory University, Atlanta, GA, United States}
\affil[c]{Merchant Ships Business Unit, Fincantieri S.p.A., Trieste, Italy}

\maketitle

\begin{abstract}
In this work, we focus on the early design phase of cruise ship hulls, where the designers are tasked with ensuring the structural resilience of the ship against extreme waves while reducing steel usage and respecting safety and manufacturing constraints.
At this stage the geometry of the ship is already finalized and the designer can choose the thickness of the primary structural elements, such as decks, bulkheads, and the shell.
Reduced order modeling and black-box optimization techniques reduce the use of expensive finite element analysis to only validate the most promising configurations, thanks to the efficient exploration of the domain of decision variables.
However, the quality of the final results heavily relies on the problem formulation, and on how the structural elements are assigned to the decision variables.
A parameterization that does not capture well the stress configuration of the model prevents the optimization procedure from achieving the most efficient allocation of the steel.
With the increased request for alternative fuels and engine technologies, the designers are often faced with unfamiliar structural behaviors and risk producing ill-suited parameterizations.

To address this issue, we enhanced a structural optimization pipeline for cruise ships developed in collaboration with Fincantieri S.p.A. with a novel data-driven hierarchical reparameterization procedure, based on the optimization of a series of sub-problems.
Moreover, we implemented a multi-objective optimization module to provide the designers with insights into the efficient trade-offs between competing quantities of interest and enhanced the single-objective Bayesian optimization module.

The new pipeline is tested on a simplified midship section and a full ship hull, comparing the automated reparameterization to a baseline model provided by the designers.
The tests show that the iterative refinement outperforms the baseline on the more complex hull, proving that the pipeline streamlines the initial design phase, and helps the designers tackle more innovative projects.
The reparameterization procedure only relies on the evaluation of surrogate models and can be applied with minimal changes to other large-scale structural problems where yielding and buckling constitute the limiting factor to the design.
\end{abstract}

\break
\tableofcontents

\section{Introduction}
\label{sec:intro}
The shipbuilding industry faces continuously evolving requirements, as the rise in environmental consciousness prescribes reducing operational costs and adopting new engine technologies~\cite{naval_green_2022}. 
The hull design should reflect a lower resource usage during the manufacturing and operational phase, while at the same time accommodating new machinery such as liquid hydrogen tanks or batteries.
All of these innovations must then observe the structural stability constraints imposed by the classification societies, which guarantee the safety and durability of the ship.
The implementation of flexible and efficient optimization procedures during the initial design phase thus represents a crucial factor for innovation and competitiveness. 
However, the formulation of the optimization problem is a delicate task. 
When facing a project with novel characteristics, the designers might not anticipate the emergence of fringe behaviors and corner cases for which their initial formulation produces underwhelming results.
In this work, we minimize the total mass of a cruise ship during the initial design phase, combining surrogates-assisted optimization and an automatic refining strategy for the optimization problem. 
The validation of the optimized designs is carried out through Finite Elements Analysis (FEA), integrating the tools and pipelines familiar to the designers.

FEA is widespread in industry~\cite{fea_industry_2008}, and commercial solvers are available for the certification of industrial designs. 
However, the computational cost of FEA is incompatible with rapid iterations of optimization steps in the initial design phase of large-scale projects.
Advanced FEA codes provide adjoint-based optimizers~\cite{opt_adjoint_2000}, but the implementation of complex and sometimes intractable rules from the classification societies hinder their adoption into the design workflow.
Recent works on the structural optimization of marine artifacts have been focused on the use of either simplified analytical formulations, or FEA of small, incomplete models.
In~\cite{naval_opt_ring_2019}, a single hull ring was optimized to reduce mass and structural instability with a multi-objective genetic algorithm, with objectives and constraints computed by FEA.
Finite differences of the simplified expressions for yielding and buckling were used to find search directions for the optimization of a semi-submersible floater in~\cite{naval_opt_floater_2019}.
FEA of a coarse mesh, combined with a simplified analytical model for the stiffeners structures, was optimized with particle swarm optimization in \cite{naval_opt_midship_2019}.
A collection of recent publications on design methods for the marine industry, including structural and topology optimization, can be found in the report of the International Ship and Offshore Structures Congress \cite{naval_opt_issc_2022}.

Reduced order models (ROMs)~\cite{rom_benner_2021, rom_certified_2016, chinestaenc2017} provide an alternative to FEA in the initial design phase, enabling the designers to quickly iterate changes to parameter values and shape configurations, with the expensive validations being performed only on the most promising configurations.
Non-intrusive approaches~\cite{rom_databook_2022, li2022machine} are able to separate the optimization procedure from the underlying physical model, a crucial requirement for workflows relying on closed-source commercial codes.
Still, the optimization results are heavily dependent on the problem formulation, which when kept simple might be hindered by the designers' bias or inexperience with novel technologies, or on the other end could scale poorly when too complex.

\begin{figure}[hbt!]
	\centering
	\includegraphics[trim=0 0 0 0, clip, width=1.\textwidth]{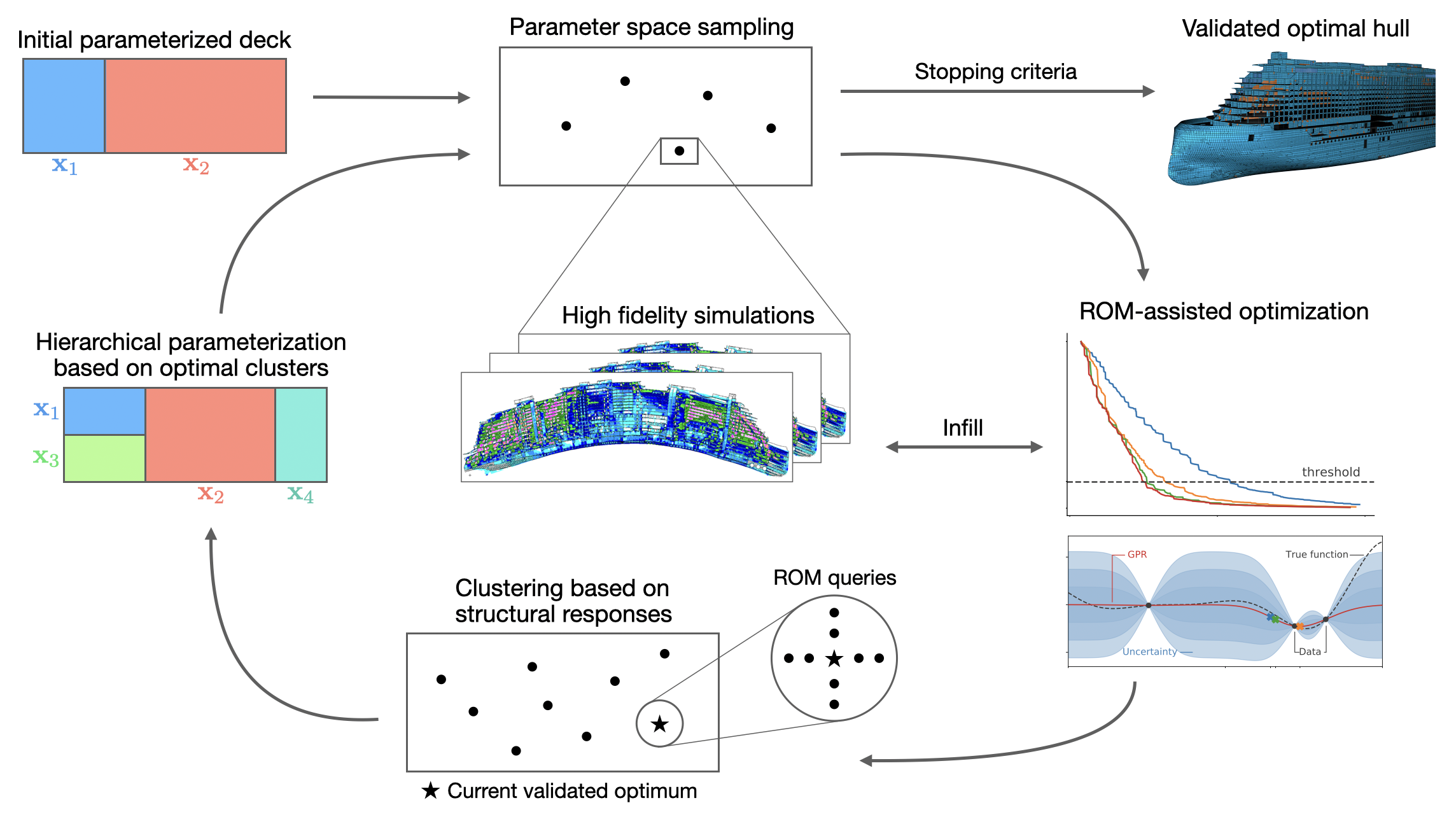}
	\caption{Optimization pipeline implementing the inner-outer loop approach.
    The process starts with a coarse parameterization of the model and a random sampling of the parametric space.
    The inner loop consists of the high-fidelity expensive simulation of the samples, from which ROMs are built in the outer loop.
    The ROMs are optimized with a multi-objective genetic algorithm and single-objective BO, with infill criteria selecting the most promising configurations for addition to the high-fidelity database.
    The reparameterization procedure starts when the optimization step is unable to find new candidates. 
    The ROMs are queried in the proximity of the current optimum to collect the structural responses, which are then clustered to generate a finer parameterization of the model, producing a new set of decision variable adapted to the emerging structural behavior.
    The new decision variables are hierarchically dependent on the previous ones and the high-fidelity database is easily updated. New samples from the larger parametric domain are drawn and the cycle repeats until the stopping criteria are met.
    }
	\label{fig:pipeline_chart}
\end{figure}

In this work, we present an automated reparameterization procedure in the context of structural optimization, where the problem formulation is refined hierarchically through an integer linear problem (ILP)~\cite{ilp_theory_1998} based on the structural responses observed near the best-known configuration.
The procedure extends an existing automated optimization pipeline for the optimization of passenger ship hulls~\cite{tezzele_multifidelity_2023}, developed in collaboration with Fincantieri S.p.A., where data-driven surrogates based on proper orthogonal decomposition (POD) and Gaussian process regression (GPR)~\cite{gpr_rasmussen_2005, gpr_gramacy_2021} are optimized using Bayesian optimization (BO)~\cite{bo_taking_2016}.
Moreover, the framework is extended with the addition of multi-objective optimization through a genetic algorithm (GA)~\cite{ga_2008, multiobj_moea_1995, multiobj_moea_2011, multiobj_nsga2_2022}, and the integration of constraints for the vertical center of gravity (VCG) in the BO procedure.
Single-objective optimization is further refined with a specialization for the discrete parametric domain, and a greedy heuristic that we call principal dimension search (PDS).
Figure~\ref{fig:pipeline_chart} shows the end-to-end pipeline, from the initial parameterization to the optimal hull. 
The process begins with an initial formulation of the optimization problem and the high-fidelity evaluation of a randomized sampling of the parametric domain.
The high-fidelity evaluations are used to build surrogates for the cheap evaluation of the quantities of interest (QoIs).
The following optimization steps, multi and single-objective, use the surrogates to efficiently select promising parameter configurations to expand the high-fidelity evaluations and refine the surrogates.
Then an automatic refinement of the problem formulation is performed, increasing the capability of reaching better values of the QoIs at the cost of introducing more decision variables.
Following the problem refinement, the sequence of surrogates construction and optimization is started again, until the user is satisfied with the QoIs of the high-fidelity optimum or some stopping criteria are met. 

The paper is organized as follows. The problem formulation for the initial design phase is described in Section~\ref{sec:problem}, along with the proposed solution for an automatic optimization pipeline.
In Section~\ref{sec:methods} we discuss the numerical methods for the implementation of the optimization pipeline.
The pipeline is applied to two test cases, a simplified midship section and a full ship model, with the results being analyzed in Section~\ref{sec:results}.
Finally, we draw conclusions in Section~\ref{sec:conclusions}.

\section{Problem definition}
\label{sec:problem}
The design phase of cruise ships can take more than one year to be completed, and is subdivided into multiple sub-phases which deal with increasing levels of detail.
In this work, we concentrate on the early design phase, where the designers are tasked with ensuring the structural resilience of the ship against extreme waves, while reducing steel usage and respecting safety and manufacturing constraints.
At this point in the project, the design of decks, shells, cabins, and functional areas has been finalized at a coarse level, defining the 2D geometry of the steel plates that constitute the primary structural components of the ship.
To modify the structural response of the model, the designers can act on the thickness and the steel type used in the primary structural components, or they can change the configuration of the secondary structural components (such as beams and flanges) by varying their size and spacing.

Every configuration needs to be validated through a computationally expensive FEA simulation.
Combined with the large variability in the design of cruise ship, which often include ad-hoc elements such as theaters, elevators disposition, or open areas, the initial design phase usually lasts several months.
In this work, we implement an optimization framework to speed up the design process and provide the designers with optimal configurations.

In Section~\ref{subsec:structural_formulation} we present the FEA setup and the QoIs extracted by a simulation. In Section~\ref{subsec:pipeline} we present the optimization pipeline.

\subsection{Structural problem formulation}
\label{subsec:structural_formulation}
\begin{figure}[hbt!]
	\centering
	\includegraphics[trim=0 0 0 0, clip, width=.45\textwidth]{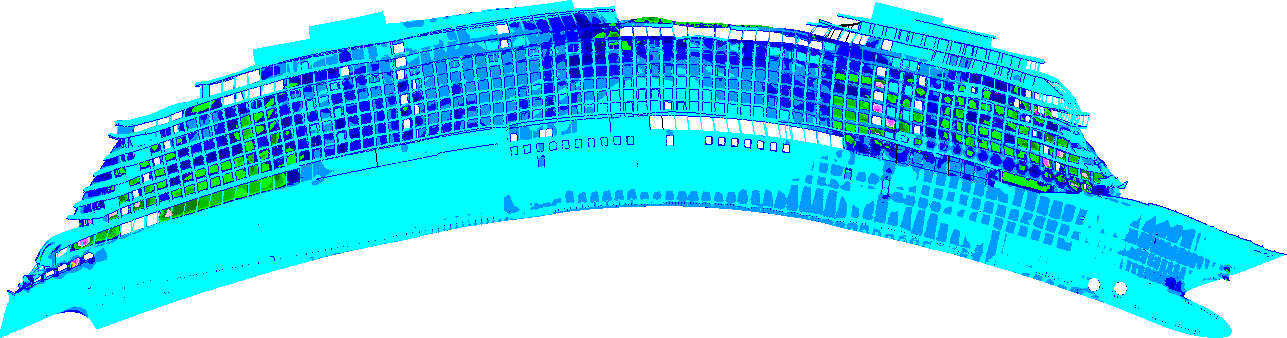}\qquad
	\includegraphics[trim=0 0 0 0, clip, width=.45\textwidth]{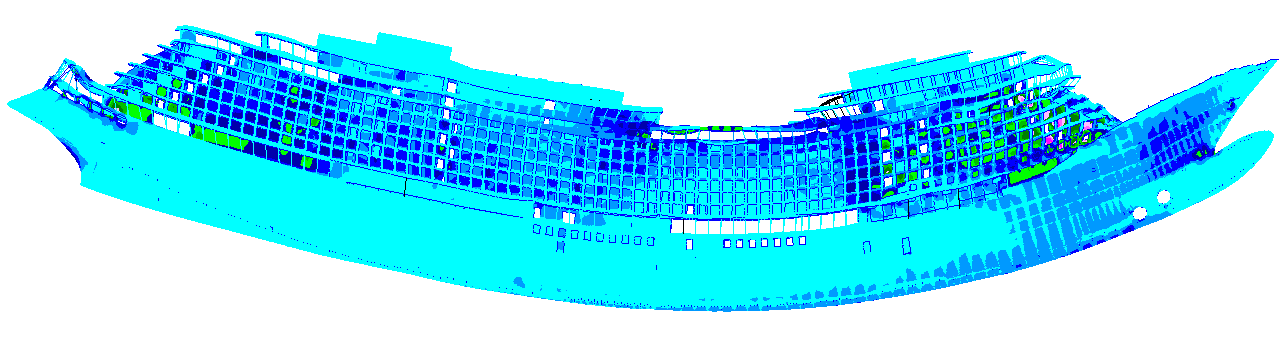}
    \caption{A full ship model under hogging load condition on the left, and sagging on the right. Displacements are magnified, colors represent the value of the von Mises yielding criterion.}
	\label{fig:patran_hogsag}
\end{figure}

The global model of a ship hull uses quadrilateral and triangular shell elements for the primary structural members (decks, bulkheads and shell), while the secondary stiffeners also use beam elements. The shell elements have principal dimensions of approximately $700$~\si{mm}.
The load conditions applied to the mesh come from two extreme waves, as defined by the classification society: the hogging condition corresponds to the wave's crest being placed at the ship's mid-length, while in the sagging condition the ship's middle point is placed above the trough. 
Figure~\ref{fig:patran_hogsag} shows the two load conditions applied to the model of a full ship. The model loads are then completed by the machinery and furniture.

MSC NASTRAN \cite{nastran_1982} assembles a linear static analysis problem for a given parameter configuration, that is a parameterized hull, using thin plate theory for the constitutive equations of shell elements.
The outputs of interest of the simulation comprise the three-dimensional displacements of the mesh nodes and the in-plane stress tensors evaluated at the centroids of the shells' faces.
The stress tensors of each shell face are averaged and transformed to the global coordinate reference system for further analysis.
At this point, the stress tensors are collected as
\begin{equation}
    S_e^l = 
    \begin{bmatrix}
        \sigma_x & \tau_{xy} & \tau_{xz} \\
        \tau_{xy} & \sigma_y & \tau_{yz} \\
        \tau_{xz} & \tau_{yz} & \sigma_z ,
    \end{bmatrix},
    \label{eq:cauchy_stress}
\end{equation}
where $l$ denotes the load condition, $e$ the shell element index, $\sigma$ the direct stresses, and $\tau$ the shear components. Subscripts to the stress components identify the axis, or plane, of action.

To quantify the structural integrity of the ship under the given load conditions, the stress tensor fields are post-processed to compute the yielding and buckling states of each element, according to the rules imposed by the classification society.
In this paper, the structural criteria come from the Det Norske Veritas (DNV) specifications for ships of length of 100 meters and above.
The yielding state is obtained by comparing the stress tensor components with the appropriate allowable value and evaluating the von Mises yielding criterion, that is:
\begin{subequations}
    \label{eq:yielding}
    \begin{align}
        |\sigma_i| \leq 245, \qquad & \text{ for } i \in \{x, y, z\}
        , \\
        |\tau_i| \leq 153, \qquad   & \text{ for } i \in \{xy, xz, yz\}
        ,
    \end{align}
    \begin{equation}
        \sigma_{\text{VM}} = \sqrt{\frac{(\sigma_x - \sigma_y)^2 + (\sigma_y - \sigma_z)^2 + (\sigma_z - \sigma_x)^2}{2} + 3(\tau_{xy}^2 +\tau_{xz}^2 + \tau_{yz}^2 )} \leq 307
        ,
    \end{equation}
\end{subequations}
where the critical values refer to high strength structural steel (AH36).
An element is marked as yielded if any of Equation~\eqref{eq:yielding} is not satisfied for at least one load condition.

Although the used shell elements do not model out-of-plane stresses, the DNV specifies buckling criteria \cite{dnv_rules} computed from the stress tensor components and the geometry of the panel the element is part of.
These criteria are based on general elastic buckling formulae, corrected in the plastic range, and result in $11$ usage factors. 
An element is marked as buckled if any of these usage factors is above a certain threshold for at least one load condition.

Both yielding and buckling phenomena can be corrected with additional manufacturing actions, albeit at the cost of increased build time and resource usage.
These corrections are undesirable, and the designers specify critical thresholds on the number of yielded and buckled elements, above which a design is to be penalized.
These thresholds provide the main constraints of the initial design phase, and the total number of yielded and buckled elements constitute the first QoIs obtained from a high-fidelity simulation.
The mass and cost of the model are computed from the volume of the steel used and the amount of adopted corrective measures. 
In this work, we only consider buckling correction through the application of secondary stiffeners, which affect both total mass and cost.

The vertical deflection is computed as the maximum absolute vertical displacement of a selected mesh node, and is interpreted as an index of the structural stiffness of the hull. 
Lower vertical deflection values correspond to more resilient ships.
The vertical center of gravity (VCG), computed for the structural components of the ship weight, is required to be lower than a critical value determined by the geometry of the hull.

The list of QoIs computed from a simulation is thus composed of five scalars: number of yielded elements, number of buckled elements, total mass considering the corrective actions, vertical deflection, and VCG.

The decision variables for the mass optimization problem are the thickness of selected shell elements.
The domain $\mathcal{T}$ for these variables is discrete, since it is limited to commercially available sizes.
A subset of the available thicknesses $\mathcal{D}_i \subseteq \mathcal{T}$ is assigned to the $i$-th group of elements to be controlled by the parameter $\xpar_i$, since regulations can specify different minimum thicknesses depending on the part of the ship. 
In this work, only primary members can be controlled by parameters, while the thickness of secondary stiffeners and the rest of the material properties are fixed.
The evaluation of a single parameters configuration requires assembling and solving the corresponding FEA problem, followed by post-processing of the results to obtain the QoIs.

\subsection{Optimization pipeline}
\label{subsec:pipeline}
We follow what has been done in \cite{tezzele_multifidelity_2023}, where the authors  
presented an automatic pipeline for the optimization of the initial ship design, using data-driven surrogates and BO to minimize the total mass of the model.
The inner-outer loop approach iterates two phases, the inner one being the high-fidelity simulation of one or more parameter configurations, and the outer phase being constituted by the selection of new configurations to be simulated, chosen through optimization of the surrogates.
The high-fidelity simulations are computationally expensive, thus it is assumed that only a limited number of runs are possible at each iteration.
The outer loop builds data-driven surrogates based on the high-fidelity simulations, with the following optimizations using a black-box approach.
The closed-source high-fidelity FEA solver used does not allow access to the implementation details, preventing intrusive approaches to model order reduction.
The loop of outer and inner operations is interrupted upon the outer phase being unable to further select promising candidates.
The end goal of our pipeline is to help the designers quickly identify the capabilities of the model, in terms of critical structural behaviors. The optimal design found by the automated procedure will necessarily undergo further analysis and modifications beyond what is modeled in our setup, but the time and resources saved in this initial phase allow the designer to only focus on the more complex problems.
In this paper, we expand the existing pipeline with additional steps in the outer phase, as depicted in Figure~\ref{fig:pipeline_chart}. We added a multi-objective optimization module and a reparameterization module. The reparameterization is the main innovation and it does not rely on the particular optimization algorithms chosen in the outer loop. 

Multi-objective optimization provides an effective tool for the exploration of the model capabilities through the discovery of optimal trade-offs between contrasting QoIs, such as total mass and vertical deflection, enabling better informed decisions on the constraints for the optimization problem.
Even for simple designs, the trade-offs between the QoIs reveal multi-modal relations between the parameters. Approaches based on random sampling of the parametric space are not able to properly identify the configurations with the optimal trade-offs, without using an overwhelming amount of samples.

The reparameterization module aims to reformulate the optimization problem itself, by increasing the parameters to improve the best-known configuration.
Indeed, the initial parameterization of the model often proves to be a limiting factor in the efficacy of the optimization, either due to being usually based on simplified beam theory or due to the model having complex or novel structures.
The reparameterization procedure identifies whether a group of elements controlled by the same parameter can be split, creating new independent parameters so that different thickness values can be assigned to the subgroups.
The following optimization is then able to reduce the mass of elements that exhibit less yielding and buckling phenomena, while increasing the stiffness of highly stressed parts of the hull.
Designers are no longer limited by the chosen initial parameterization, where a modification would require a tedious manual reconfiguration and data transfer. 

Finally, a number of improvements to the single-objective optimization module enable a more efficient exploration of the parametric domain, which becomes necessary to handle the growing number of parameters coming from the reparameterization step.

\section{Numerical methods}
\label{sec:methods}
This section presents the numerical methods used in the automatic optimization pipeline.
Section~\ref{subsec:rom} discusses the implementation of data-driven ROMs to compute the QoIs, Section~\ref{subsec:opt} illustrates the optimization procedures enabled by the surrogates, and Section~\ref{subsec:reparam} describes the parameterization refinement process.

\subsection{Model order reduction}
\label{subsec:rom}
To implement a computationally cheap and accurate surrogate of the stress tensor field for different load conditions, we combine POD and GPR.
In this subsection, we briefly present these two methods.

\subsubsection{Proper orthogonal decomposition}
\label{subsubsec:pod}
POD decomposes a matrix of high-fidelity snapshots to obtain an orthogonal basis for the span of the column space.
Let $\mathbf{M} \in \mathbb{R}^{n \times m}$ be a matrix collecting the $m$ snapshots as columns of $n$ features, and assume $n \gg m$ which corresponds to operating in a large-scale scenario.
Then, the singular value decomposition (SVD) of $\mathbf{M}$, that is $\mathbf{M} = \mathbf{U} \mathbf{\Sigma} \mathbf{V}^T$,
uniquely determines $\mathbf{U} \in \mathbb{R}^{n \times n}$ the matrix with the left singular vectors of $\mathbf{M}$ as columns, $\mathbf{\Sigma} \in \mathbb{R}^{n \times m}$ the diagonal matrix of singular values, and $\mathbf{V}^T \in \mathbb{R}^{m \times m}$ the transpose of the matrix with the right singular vectors of $\mathbf{M}$ as columns.
Uniqueness is achieved by sorting the singular values in decreasing order, up to a change of sign of the basis vectors or a change in positions corresponding to repeated singular values.
By choosing a truncation rank $r < m$, one can introduce the truncated SVD 
\begin{equation}
    \label{eq:svd_truncated}
    \mathbf{M} \approx \tilde{\mathbf{M}} = \tilde{\mathbf{U}} \tilde{\mathbf{\Sigma}} \tilde{\mathbf{V}}^T
    ,
\end{equation}
which constructs $\tilde{\mathbf{U}}$ and $\tilde{\mathbf{V}}$ by selecting the $r$ leftmost columns of $\mathbf{U}$ and $\mathbf{V}$, and the corresponding $r$ largest singular values form the diagonal matrix $\tilde{\mathbf{\Sigma}}$.
The reconstruction error, computed as the Frobenius norm of $\mathbf{M} - \tilde{\mathbf{M}}$, is given by the sum of the squares of the $m-r$ discarded singular values.
The columns of $\tilde{\mathbf{U}}$ thus give a basis of size $r$, optimal in the sense of the Frobenius norm, for the reconstruction of the original snapshots.
The reduced coefficients can be computed through the projection of the snapshots matrix $\mathbf{M}$ unto the basis $\tilde{\mathbf{U}}$, with the matrix $\mathbf{C} = \tilde{\mathbf{U}}^T \mathbf{M}$ holding as $i$-th column the reduced coefficients of for the $i$-th snapshot in $\mathbf{M}$.
As each snapshot is associated with a parameter configuration, the reduced coefficients provide the ground truth for the construction of a map from the parametric domain to the space of reduced coefficients.
In this paper, the snapshots are obtained through the execution of closed-source commercial codes, thus regression can only be carried out by data-driven methods.
Successful approaches that have been leveraged in the past are radial basis functions interpolation, Gaussian process regression and artificial neural networks, among others~\cite{swischuk2019projection, demo2021hull, pod_dl_rom_2022, gpr_pod_guo_2018}.

\subsubsection{Gaussian process regression}
\label{subsubsec:gpr}
GPR assumes that any set of observations of the QoI is sampled from a multivariate normal distribution.
Let $\ypar$ be a vector of observations $\ypar_i = f(\xpar^{(i)})$, and $\ypar \sim \mathcal{N} \left( \mupar, \mathbf{K} \right)$.
The vector of marginal means $\mupar$ and the covariance matrix $\mathbf{K}$ are defined as:
\begin{gather}
    \mupar_i = \mathbb{E}\left[ \ypar_i \right] = \mathbb{E}\left[ f(\xpar^{(i)}) \right]
    , \\
	\mathbf{K}_{ij} = \mathbf{K}_{ji} = \mathbb{E}\left[ (\ypar_i - \mupar_i)(\ypar_j - \mupar_j) \right] = \mathbb{C}\text{ov}\left[ \ypar_i, \ypar_j \right] = \mathbb{C}\text{ov}\left[ f(\xpar^{(i)}), f(\xpar^{(j)}) \right]
    .
\end{gather}
The conditional distribution of a subset of observations $L$, on past outcomes $H$, can be obtained by leveraging the joint distribution
\begin{equation}
    \label{eq:gp_conditional_full}
	\begin{bmatrix*}[l]
		\ypar_H \\
		\ypar_L \\
	\end{bmatrix*} 
	\sim \mathcal{N} \left( 
		\begin{bmatrix*}[l]
			\mupar_H \\
			\mupar_L \\ 
		\end{bmatrix*}
		,
		\begin{bmatrix*}[l]
			\mathbf{K}_{HH} & \mathbf{K}_{HL} \\
			\mathbf{K}_{LH} & \mathbf{K}_{LL} \\
		\end{bmatrix*}
	\right)
    ,
\end{equation}
to obtain 
\begin{equation}
    \label{eq:gp_y_conditioned}
	\ypar_{L \mid H} \sim \mathcal{N}(\mupar_{L \mid H}, \mathbf{K}_{L \mid H})
    ,
\end{equation}
where the conditioned terms are expressed in terms of mean and covariance matrix of the process' distribution as
\begin{equation}
    \label{eq:gp_mean_conditioned}
	\mupar_{L \mid H} = \mupar_L + \mathbf{K}_{LH} \mathbf{K}_{HH}^{-1} (\ypar_H - \mupar_H)
    ,
\end{equation}
\begin{equation}
    \label{eq:gp_cov_conditioned}
	\mathbf{K}_{L \mid H} = \mathbf{K}_{LL} - \mathbf{K}_{LH} \mathbf{K}_{HH}^{-1} \mathbf{K}_{HL}
    .
\end{equation}

The crucial component of a GPR is the construction of two functions $\mu(\xpar^{(i)})$ and $\text{kern}(\xpar^{(i)}, \xpar^{(j)})$ to generate the entries of $\mupar$ and $\mathbf{K}$, respectively, so that past observations are reproduced with high confidence and extrapolation beyond the observed instances gives a coherent probability distribution.
The mean function, which encodes any prior knowledge about the output, is usually set to zero.
Indeed, in Equation~\eqref{eq:gp_mean_conditioned} the covariance-dependent terms act as a correction, weighting the observed deviations from the prior.
For the covariance function, which is required to generate a symmetric semi-definite positive matrix, a popular choice is the squared exponential kernel with automatic relevance determination:
\begin{equation}
    \label{eq:gp_cov_exp}
    \mathbf{K}_{ij} = \text{kern}(\xpar^{(i)}, \xpar^{(j)}) = \sigma^2 \text{exp}\left( -\frac{1}{2}\sum_{d}\frac{(\xpar^{(i)}_d - \xpar^{(j)}_d )^2}{l_d^2} \right)
    ,
\end{equation}
where $\sigma^2$ is a scaling factor, and $l_d$ is the length scale associated to the dimension $d$ of the parametric domain.

To fit the covariance matrix $\mathbf{K}_{\bm\theta}$ parameterized by the vector of hyperparameters $\bm\theta$, the log likelihood is maximized through gradient ascent.
For the case of zero prior and the kernel function in Equation~\eqref{eq:gp_cov_exp}, with $\Xpar = \left[ \xpar^{(1)}, \dots , \xpar^{(m)} \right]$, it reads
\begin{equation}
    \label{eq:gp_loglikelihood}
    \log p(\ypar \mid \Xpar, \bm\theta ) = -\frac{m}{2}\log 2\pi - \frac{1}{2} \log | \mathbf{K}_{\bm\theta} | -\frac{1}{2} \ypar^T  \mathbf{K}_{\bm\theta}^{-1} \ypar
    ,
\end{equation}
where $| \mathbf{K}_{\bm\theta} |$ is the determinant of the covariance matrix.

We emphasize that it is possible to precompute $\mathbf{K}_{HH}^{-1} (\ypar_H - \mupar_H)$ in Equation~\eqref{eq:gp_mean_conditioned}, since it only depends on the snapshots.
Thus, the online query for the conditioned mean requires only the computation of the kernel function between the snapshots samples, and the query point.
The execution time of a GPR query on a single point scales with the product of the number of snapshots samples, output dimensions, and kernel function complexity.

\subsubsection{Implementation of the surrogates}
\label{subsubsec:surrogates}
In this paper, our aim is to find a computationally inexpensive map from the parameters configuration to the stress tensor fields.
We leverage POD to reduce the dimensionality of the output space, and use GPR to map the parameter configurations to the vector of reduced coefficients, thus obtaining a much simpler problem setup.
This approach is known as POD-GPR in the literature, and has been successfully used in solid mechanics \cite{gpr_pod_guo_2018, pod_gpr_crash_2022, pod_gpr_uq_2023}.
In the following, we omit the indices for the load condition and the stress tensor component to simplify the notation.

Let $\mathbf{S} \in \mathbb{R}^{n \times m}$ be the matrix collecting the snapshots vectors $\mathbf{s}^{(i)} := \mathbf{s} (\mathbf{x}^{(i)}) $ of the current stress tensor component and load condition, associated with the parameter configuration $\xpar^{(i)}$.
Let $\tilde{\mathbf{U}}$ be the basis matrix obtained by the SVD of $\mathbf{S}$, truncated at rank $r$.
We obtain the parameterized reduced coefficients through the projection 
\begin{equation}
    \mathbf{c}^{(i)} := \mathbf{c} (\mathbf{x}^{(i)}) = \tilde{\mathbf{U}}^T \mathbf{s}^{(i)}
    ,
\end{equation}
For each vector of $r$ reduced coefficients, we train a different GPR on the $m$ pairs $(\xpar^{(i)}, \mathbf{c}^{(i)})$, so that the approximated quantities, denoted by a hat, for the coefficients and the original field satisfy
\begin{equation}
    \mathbf{s}^{(i)} \approx \hat{\mathbf{s}}^{(i)} = \tilde{\mathbf{U}} \hat{\mathbf{c}}^{(i)}
    ,
\end{equation}
where the hat denotes the result of an approximation. 

In this paper, we consider 2 load conditions (hogging and sagging) and the unique 6 elements of the Cauchy stress tensor Equation~\eqref{eq:cauchy_stress}, therefore the total number of vector-valued GPRs being fitted is 12.
We use vector-valued GPRs, each mapping a parameter configuration to the $r$ POD coefficients corresponding to a different combination of load condition and stress tensor component.
We also observe that the decay of singular values consistently achieves residual energy below $1\%$ for rank $r$ higher than, but close to, the number of parameters.
This behavior results in an effective model order reduction and results in a cheap training phase for the GPRs.
Additionally, the construction of the GPRs for each stress tensor component and load condition is embarrassingly parallel, once the tensor fields are decomposed and each component stored separately.

The post-processing steps for the computation of the yielding and buckling states of each element only require the stress tensor components and metadata of the element, leading to an efficient vectorized implementation.
The same holds for the aggregation to compute the total number of yielded and buckled elements, from which the derivation of total mass and cost is straightforward.

We remark that the use of an efficient post-processing step, instead of the direct implementation of surrogates for the QoIs, results in a much lower error in the prediction, as was shown in~\cite{tezzele_multifidelity_2023}.
For the surrogate of the vertical deflection, we construct a single GPR for each load condition and, during a query, take the maximum of the absolute value of the predictions.

The vector of QoIs computed from a parameter configuration $\xpar$ is thus
\begin{equation}
\label{eq:qois}
[n_\text{y}(\xpar), n_\text{b}(\xpar), f_\text{deflection}(\xpar), f_\text{mass}(\xpar), \text{VCG}(\xpar)]
\end{equation}
where $n_\text{y}(\xpar)$ is the number of elements for which the predicted stress tensor, under any load condition, fails at least one of Equations~\eqref{eq:yielding}. The number of buckled elements $n_\text{b}(\xpar)$ follows the same principle, but comes from the aggregation of the buckling usage factors computed according to~\cite{dnv_rules}. $f_\text{deflection}(\xpar)$ is the vertical deflection, while the mass is computed as
\begin{equation}
    \label{eq:mass}
    f_\text{mass}(\xpar) = m_{\text{fixed}} + \mathbf{d} \cdot \xpar + m_{\text{bar}} n_\text{b}(\xpar)
    ,
\end{equation}
where $m_{\text{fixed}}$ represents the mass of the elements not controlled by any parameter, $\mathbf{d}$ is the vector collecting the linear density of each parameterized section, and $m_{\text{bar}}$ is the average mass of a reinforcement bar.
Finally, the VCG is computed as 
\begin{equation}
    \label{eq:vcg}
    \text{VCG}(\xpar) = \frac{\text{VCG}_{\text{fixed}} m_{\text{fixed}} + \sum_{i} \text{VCG}_{i} \mathbf{d}_i \xpar_i}{m_{\text{fixed}} + \sum_{i} \mathbf{d}_i \xpar_i} \geq 0
    ,
\end{equation}
where $\text{VCG}_{\text{fixed}}$ is the VCG of elements not controlled by any parameter, and $\text{VCG}_{i}$ is the VCG of the $i$-th parameterized region. 
VCG is conventionally assumed to be positive.

Overall, the computational complexity of querying our surrogates scales as $\mathcal{O}(r(n+md))$, where $d$ is the number of parameters.
With both $m \ll n$ and $d \ll n$, the cost of a surrogate query is much lower than the $\mathcal{O}(n^2)$ required by a high-fidelity FEA for a sparse linear system.
This enables the use of the surrogates in many-query and quasi-real time applications, such as black-box optimization procedures and graphical user interfaces.

\subsection{Surrogate-assisted optimization}
\label{subsec:opt}
The surrogates for the stress field, combining POD and GPR, allow a large number of inexpensive queries on the parametric domain. 
Only the most favorable configurations, selected by optimization, will be validated with an expensive full order simulation. 
However, as in many industrial applications, the parametric domain is large enough that exhaustive exploration is still unfeasible, and an infill criterion is required. 
We propose a combination of a multi-objective genetic algorithm, Bayesian optimization, and a heuristic local search to guide the exploration of the parametric domain and perform full order simulations only for the most promising configurations.

In Section~\ref{subsubsec:multiobj_opt}, we present the genetic algorithm for multi-objective optimization, and an infill criterion which leverages the uncertainty quantification provided by the GPR component of our surrogates.
Section~\ref{subsubsec:opt_bayesian} presents a specialization of BO, discussing the choice of objective function and the implementation of additional constraints and heuristics.
Section~\ref{subsubsec:opt_pds} illustrates the black-box local heuristic used to further refine the optimum found by BO.

\subsubsection{Multi-objective optimization}
\label{subsubsec:multiobj_opt}
Often, multiple QoIs are used as objective functions and it is not clear whether a global criterion for the prioritization of one over another exists.
In this case, the experts are interested in obtaining not a single optimal configuration, but rather a set of configurations exhibiting the best trade-offs between the different QoIs. 
Such a set is called the Pareto set (PS) of the problem, and its image in the space of QoIs is the Pareto frontier (PF). 
By analyzing the PF, the experts can compromise between the various needs and impose stricter or laxer constraints for the further single-objective optimization tasks. 

Let \Dpar~be a parametric domain, and $f$ a vector of QoIs $\{f_i: \Dpar \rightarrow \mathbb{R} \}_{i=1}^n$ which must be all minimized. 
A configuration $\bar{\xpar} \in \Dpar$ is called dominated if there is at least one $\xpar \in \Dpar$ such that $f_i(\xpar) \leq f_i(\bar{\xpar})$ for all $i$, and $f_i(\xpar) < f_i(\bar{\xpar})$ for at least one $i$. 
A dominated configuration is of no interest in the optimization since there exists another configuration that performs no worse on all QoIs, and is strictly better for at least one QoI.
Given a population of configurations, the PS is then taken as the subset of non-dominated individuals.

In our framework, we use multi-objective optimization of the QoIs shown in Equation~\eqref{eq:qois} to enrich the initial sampling of high-fidelity simulations and provide the designers with an approximated PF.
The designers can then analyze the PF and decide on the critical thresholds for the following single-objective optimization, by estimating the heuristic trade-offs between mass and structural QoIs.

To approximate the PF, many popular approaches use genetic algorithms \cite{ga_2008, multiobj_moea_1995, multiobj_moea_2011} in which an initial population is iteratively grown through crossover and mutation of its best individuals, and the least desirable elements are selected for removal.
Other approaches from the literature include the construction of the PS by collaboration between multiple random scalarizations of the QoIs \cite{multiobj_moead_2007}.
A number of methods have been developed for probabilistic surrogates, extending single-objective Bayesian Optimization to the multi-objective setting \cite{mobo_parego_2006, mobo_scalarize_2021, mobo_morbo_2022}.
These methods usually maintain multiple GPRs built from a sample population, where the addition of new samples comes by optimizing some measure of quality of the PF.
The usual choice is maximizing the dominated hyper-volume, as with the Expected HyperVolume Improvement (EHVI) \cite{mobo_ehvi_2006}.
However, these methods add a further level of surrogation, require Monte Carlo evaluation of the EHVI against the current PF and the acquisition procedure needs to be adapted to the discrete parametric domain of our problem.
We chose to adopt a genetic algorithm to directly query the surrogates described in Section~\ref{subsubsec:surrogates}, which support efficient batch evaluations.
The parametric domain is explored by modifying the population according to the crossover, mutation and selection operators, which can are implemented to guarantee the feasibility and diversity of the samples.
The generic structure is outlined in Algorithm~\ref{alg:genetic}.

\begin{algorithm}[ht]
\begin{algorithmic}
    \Require initial population $\Xpar$, fitness function $f$, crossover and mutation functions, maximum population $p$, number of generations $q$
    \Ensure final population
    \State $i \gets 0$
    \While{$i < q$}
        \State $\Xpar^\text{parents} \gets$ select from $\Xpar$ according to $f$
        \State $\Xpar^\text{children} \gets$ crossover of $\Xpar^\text{parents}$
        \State $\Xpar^\text{children} \gets$ mutation of $\Xpar^\text{children}$ 
        \State $\Xpar \gets \Xpar^\text{parents} \cup \Xpar^\text{children}$
        \State $\Xpar \gets$ the $p$ best performing individuals from $\Xpar$ according to $f$
        \State $i \gets i + 1$
    \EndWhile \\
    \Return $\Xpar$
\end{algorithmic}
    \caption{Template of a genetic algorithm.}
    \label{alg:genetic}
\end{algorithm}

For the selection of fittest individuals, we use the non-dominated sorting strategy proposed in the NSGA family of algorithms \cite{multiobj_nsga2_2022, multiobj_nsga3_2014}.
Individuals are ranked by their degree of domination with Algorithm~\ref{alg:ns-sort}.

\begin{algorithm}
\begin{algorithmic}
    \Require list of objective functions evaluations $\Ypar = \left\{\ypar^{(i)} \right\}_{i=1}^m$
    \Ensure list $L$ of non-dominated layers as set of indices
    \State $L \gets $ empty list
    \State $R \gets \{ i=1, ..., m\} $
    \While{$|R| > 0$}
        \State $\Ypar^\text{R} \gets \{ \ypar^{(i)} \mid i \in R \}$
        \State $M \gets \{ i \mid \ypar^{(i)} \text{  is not dominated  in } \Ypar^\text{R} \}$
        \State $L \gets L \cup M$ 
        \State $R \gets R \setminus M$
    \EndWhile \\
    \Return $L$
\end{algorithmic}
    \caption{Non-dominated sorting.}
    \label{alg:ns-sort}
\end{algorithm}

The selection of the fittest individuals is performed by taking the union of non-dominated subsets, starting from the first, until the limit of population size is reached.
If the number of selected individuals exceeds the limit, the subset with the most dominated individuals is downsampled with the niching strategy from NSGA-III~\cite{multiobj_nsga3_2014}, to prevent overcrowding in the QoIs space.

Due to the use of surrogates, the final PF is only an approximation of the true set of non-dominated QoIs.
An infill criterion is needed to select a reduced number of individuals to be validated by the high-fidelity solver, so that further multi-objective optimizations yield a more accurate PF.
The approximated PS comprises the high-fidelity samples already validated and which contribute to the construction of the surrogates, and the low-fidelity samples which can be added to the high-fidelity set to increase the accuracy of the surrogates. 
For an individual to be chosen, it is required that its addition to the high-fidelity samples maximally improve the predictions at the location of the remaining low-fidelity members of the PF.
This approach is similar to \cite{infill_alc_2000}, which adapts the idea of maximizing the conditioned uncertainty reduction from~\cite{infill_alm_1992} and the integration of an uncertainty measure over a domain of interest from \cite{infill_cohn_1996}, to the case of GPR.
This principle has been developed in terms of entropy in \cite{mobo_infill_pesmo_2016, mobo_infill_mesmo_2019, mobo_infill_pfes_2020} for multi-objective optimization for GPRs that directly model the objective functions.
We choose to measure the covariance between the samples in the reduced coefficient space since the corresponding GPRs have been thoroughly validated as providing accurate predictions of the stress tensors, from which the QoIs are derived.
Differently from the procedure in \cite{infill_alc_2000}, which would explicitly estimate the reduction in uncertainty, we estimate the increment in covariance, which only requires the evaluation of the kernel and is thus cheap compared with Equation~\eqref{eq:gp_cov_conditioned}.
Moreover, we operate on a set of GPRs so that an aggregation strategy becomes necessary.

The infill procedure starts by constructing the covariance matrices between the elements in $\Xpar^\text{L}$, and between those in $\Xpar^\text{L}$ and $\Xpar^\text{H}$. 
We aggregate over reduced coefficients and load conditions by taking the maximum, then sum over the stress tensor components to obtain $\mathbf{C}^\text{LL} \in \mathbb{R}^{n \times n}$, symmetric, and $\mathbf{C}^\text{LH} \in \mathbb{R}^{n \times m}$.

Let $\{ \mathbf{K}_\text{LL}^{(l, s, c)} \}$ be the set of covariance matrices of the GPR corresponding to the reduced coefficient $c$ of the stress tensor component $s$ under the load condition $l$, evaluated between the samples in $\Xpar^\text{L}$.
Then, the element $(i, j)$ of the aggregation matrix is computed as
\begin{equation}
    \label{eq:mobj_infill_aggregation}
    \mathbf{C}^\text{LL}_{ij} = \sum_{c} \max_{l, s} \left [ \left(\mathbf{K}_\text{LL}^{(l, s, c)}\right)_{ij} \right ],
    \qquad \forall \, i, j = 1, \ldots,  n
    ,
\end{equation}
giving a scalar measure of the mutual information between the samples in $\Xpar^\text{L}$.
The same expression is used for the construction of $\mathbf{C}^\text{LH}$.
A vector $\Delta$ is built by computing, for the $i$-th sample in $\Xpar^\text{L}$, the total positive relative increase over the maximum high-fidelity contribution to the rest of the low-fidelity PF as
\begin{equation}
    \Delta_i = \frac{\sum\limits_{j \not= i} ( \mathbf{C}^\text{LL}_{ij} - \max\limits_{h} \mathbf{C}^\text{LH}_{jh} )_+}{\sum\limits_{j \not= i} \max\limits_{h} \mathbf{C}^\text{LH}_{jh}},
    \qquad \forall \, i = 1, \ldots,  n
    ,
    \label{eq:mobj_infill_delta}
\end{equation}
where $( \cdot )_+ = \max(\cdot, 0)$.
The high-fidelity set is then enriched by adding the $i^*$-th sample determined~by
\begin{equation}
    \label{eq:mobj_infill_argmax}
    i^* = \argmax_{i} \Delta
    .
\end{equation}

The choice of the aggregation functions in Equation~\eqref{eq:mobj_infill_aggregation} can be viewed as constructing the best case scenario for the mean covariance increase across the stress tensor components.
If more than one configuration can be selected for the high-fidelity validation, the matrices $\mathbf{C}^\text{LL}$ and $\mathbf{C}^\text{LH}$ can be updated by simulating the addition of the sample $i^*$ to the high-fidelity set, as summarized in Figure~\ref{fig:mobj_infill_illustration}. 
In particular, the $i^*$-th row (and column) is removed from $\mathbf{C}^\text{LL}$, its elements are used to build the column $m+1$ in $\mathbf{C}^\text{LH}$, and finally the $i^*$-th row of $\mathbf{C}^\text{LH}$ is removed. 
The resulting matrices have thus size $n-1 \times n-1$ and $n-1 \times m+1$ respectively, and the selection can take place after computing the updated $\Delta$ vector.

\pgfmathdeclarerandomlist{MyRandomColors}{%
    {blue!20!cyan!10}%
    {blue!40!cyan!60}%
    {blue!40!cyan!50}%
    {blue!40!cyan!40}%
    {blue!30!cyan!60}%
    {blue!30!cyan!50}%
    {blue!30!cyan!40}%
    {blue!30!cyan!30}%
    {blue!20!cyan!40}%
    {blue!20!cyan!30}%
    {cyan!20}%
}%

\newcommand*{\ColorCells}[3]{%
    \pgfmathsetseed{#3}
    \foreach \y in {1,...,#1} {
        \foreach \x in {1,...,#2} {
            \pgfmathrandomitem{\RandomColor}{MyRandomColors} 
            \draw [fill=\RandomColor, fill opacity=0.7, line width=0.11mm, scale=0.33] 
                (\x-1, \y-1) rectangle (\x, \y);
        }%
    }%
}%

\newcommand*\ColorCellsSymmetric[3]{
    \def\nrows{#2}
    \def\ncols{#1}

    \foreach \i in {0,...,\numexpr \nrows-1 \relax} {
        \pgfmathsetmacro{\y}{#1 -\i}
        \pgfmathsetmacro{\x}{\i +1}
        
        \draw [fill=blue!60!cyan!60, fill opacity=0.7, line width=0.11mm, scale=0.33] 
            (\x-1, \y-1) rectangle (\x, \y);
    }%
    
    \pgfmathsetseed{#3}
    \foreach \i in {0,...,\numexpr \nrows-2 \relax} {
        \foreach \j in {\numexpr \i+1 \relax,...,\numexpr \ncols-1 \relax} {
            \pgfmathrandomitem{\RandomColor}{MyRandomColors} 

            \pgfmathsetmacro{\y}{\ncols -\i}
            \pgfmathsetmacro{\x}{\j +1}
            \draw [fill=\RandomColor, fill opacity=0.7, line width=0.11mm, scale=0.33] 
                (\x-1, \y-1) rectangle (\x, \y);

            \pgfmathsetmacro{\y}{\ncols -\j}
            \pgfmathsetmacro{\x}{\i +1}
            \draw [fill=\RandomColor, fill opacity=0.7, line width=0.11mm, scale=0.33] 
                (\x-1, \y-1) rectangle (\x, \y);
        }%
    }%
}%

\begin{figure}[!htb]
\centering
\begin{tikzpicture}[every node/.style={outer sep=0pt}]
    
    \begin{scope}[thick, xshift=1cm] 
        \ColorCellsSymmetric{10}{10}{42}
        
        \draw[scale=0.33, color=blue!90!cyan!70, line width=.5mm] (0,3) rectangle (6,4);
        \draw[scale=0.33, color=blue!90!cyan!70, line width=.5mm] (7,3) rectangle (10,4);

        \node at (1.65, 3.65) {$\mathbf{C}^{\text{LL}}$};
        \node at (-0.25, 1.16) {\footnotesize $i^*$};
    \end{scope}

    \begin{scope}[thick, xshift=5cm] 
        \ColorCells{10}{3}{42}
        
        \draw[scale=0.33, color=blue!90!cyan!70, line width=.5mm] (2,0) grid (3,1);
        \draw[scale=0.33, color=blue!90!cyan!70, line width=.5mm] (1,1) grid (2,2);
        \draw[scale=0.33, color=blue!90!cyan!70, line width=.5mm] (0,2) grid (1,3);
        \draw[scale=0.33, color=blue!90!cyan!70, line width=.5mm] (2,4) grid (3,5);
        \draw[scale=0.33, color=blue!90!cyan!70, line width=.5mm] (0,5) grid (1,6);
        \draw[scale=0.33, color=blue!90!cyan!70, line width=.5mm] (0,6) grid (1,7);
        \draw[scale=0.33, color=blue!90!cyan!70, line width=.5mm] (2,7) grid (3,8);
        \draw[scale=0.33, color=blue!90!cyan!70, line width=.5mm] (1,8) grid (2,9);
        \draw[scale=0.33, color=blue!90!cyan!70, line width=.5mm] (2,9) grid (3,10);

        \node at (.5, 3.65) {$\mathbf{C}^{\text{LH}}$};
    \end{scope}

    \begin{scope}[thick, xshift=6.66cm] 
        \ColorCells{10}{1}{42}
        
        \draw[scale=0.33, color=blue!90!cyan!70, line width=.5mm] (0,3) grid (1,4);

        \node at (.16, 3.6) {$\Delta$};
    \end{scope}

    \node [text width=1.8cm, align=center] at (7.8, 1.16) {\footnotesize max infill criterion};
    \draw [-stealth] (7.4, 2.4) -- (9., 2.4);
    \node [text width=1.8cm, align=center] at (8.2, 2.66) {\small new HF};

    \begin{scope}[thick, xshift=9.5cm] 
        \ColorCellsSymmetric{10}{10}{42}
        
        \draw[scale=0.33, color=purple!70, line width=.5mm] (0,3.5) -- (10,3.5);
        \draw[scale=0.33, color=purple!70, line width=.5mm] (6.5,0) -- (6.5,10);

        \node at (1.65, 3.65) {$\mathbf{C}^{\text{LL}}$};
        \node at (-0.25, 1.16) {\footnotesize $i^*$};
    \end{scope}

    \begin{scope}[thick, xshift=13.5cm] 
        \ColorCells{10}{3}{42}
        \draw [fill=blue!30!cyan!30, fill opacity=0.7, line width=0.11mm, scale=0.33] (3, 9) rectangle (4, 10);
        \draw [fill=blue!20!cyan!10, fill opacity=0.7, line width=0.11mm, scale=0.33] (3, 8) rectangle (4, 9);
        \draw [fill=cyan!20, fill opacity=0.7, line width=0.11mm, scale=0.33] (3, 7) rectangle (4, 8);
        \draw [fill=blue!20!cyan!10, fill opacity=0.7, line width=0.11mm, scale=0.33] (3, 6) rectangle (4, 7);
        \draw [fill=cyan!20, fill opacity=0.7, line width=0.11mm, scale=0.33] (3, 5) rectangle (4, 6);
        \draw [fill=blue!30!cyan!30, fill opacity=0.7, line width=0.11mm, scale=0.33] (3, 4) rectangle (4, 5);
        \draw [fill=blue!60!cyan!60, fill opacity=0.7, line width=0.11mm, scale=0.33] (3, 3) rectangle (4, 4);
        \draw [fill=blue!30!cyan!40, fill opacity=0.7, line width=0.11mm, scale=0.33] (3, 2) rectangle (4, 3);
        \draw [fill=blue!30!cyan!60, fill opacity=0.7, line width=0.11mm, scale=0.33] (3, 1) rectangle (4, 2);
        \draw [fill=cyan!20, fill opacity=0.7, line width=0.11mm, scale=0.33] (3, 0) rectangle (4, 1);

        \draw[scale=0.33, color=blue!90!cyan!70, line width=.5mm] (3,0) rectangle (4,10);
        
        \draw[scale=0.33, color=purple!70, line width=.5mm] (0,3.5) -- (4,3.5);

        \node at (0.75, 3.65) {$\mathbf{C}^{\text{LH}}$};
    \end{scope}

  \end{tikzpicture} 

    \caption{On the left, position $i^*$ in the low-fidelity Pareto Frontier is selected as maximizer of the infill criterion $\Delta$ evaluated on matrices $\mathbf{C}^\text{LL}$ and $\mathbf{C}^\text{LH}$ using Equation~\eqref{eq:mobj_infill_delta}. The addition of the $i^*$-th sample to the high-fidelity set is simulated with the removal of row (and column) $i^*$ from $\mathbf{C}^\text{LL}$, its addition to $\mathbf{C}^\text{LH}$, followed by the removal of row $i^*$. The next sample is selected using the updated matrices.}
    \label{fig:mobj_infill_illustration}
\end{figure}

Assuming that the covariance function $\text{kern}(\cdot, \cdot)$ of the GPRs remains the same after fitting the updated high-fidelity snapshots, the covariance matrix blocks $\mathbf{K}_{\text{HH}}$ and $\mathbf{K}_{\text{HL}}$ that appear in Equation~\eqref{eq:gp_conditional_full} will contain larger entries in the rows and columns corresponding to the elements from the infill set.
Thus, when queried on $\xpar$ from the surrogate PS, a more effective correction of the prior will be applied in Equation~\eqref{eq:gp_mean_conditioned} and likewise a larger reduction of the variance will come from Equation~\eqref{eq:gp_cov_conditioned}.
In practice, this infill criterion proves effective in generating a sparse sampling of the surrogate PS, and the rate of change of $\Delta_{i^*}$ provides a convergence criterion for early stopping of the multi-objective optimization.

\subsubsection{Bayesian optimization}
\label{subsubsec:opt_bayesian}
BO aims to optimize an expensive black box objective function, using a low number of function evaluations.
The procedure builds a surrogate and updates it iteratively, selecting the next sample as the most promising for improving on the best-known configuration. 
A GPR is built on an initial set of objective function evaluations, and the selection of the samples is carried out by optimizing an easy-to-compute acquisition function, which leverages both the predicted value and the associated uncertainty \cite{gpr_gramacy_2021}.
The sampling of the parametric domain refines the GPR, thus effectively finding good configurations. 
This sampling balances between the so-called exploration (regions with high uncertainty) and exploitation (regions close to the current optimum).

The acquisition function represents the uncertain gain in terms of objective function, so its maximization will select the next sample for the validation.
Several formulations have been proposed, which vary in terms of exploratory or exploitative proclivity.
Here, we briefly describe the negative lower confidence bound, the expected improvement, and the probability of improvement.

Negative lower confidence bound (NLCB) \cite{bo_nlcb_2009} is an exploratory acquisition function for minimization problems, expressed as
\begin{equation}
    \label{eq:acq_nlcb}
    \alpha_{\text{NLCB}}(\xpar) = -(\mu(\xpar) - \beta\sigma(\xpar))
    ,
\end{equation}
where $\mu$ is the mean and $\sigma$ is the standard deviation from the GPR posterior, with $\beta \geq 0$ weighting optimistically the contribution from the estimated uncertainty. 

Expected improvement (EI) \cite{bo_ei_2020} is an exploratory acquisition function defined as
\begin{equation}
    \label{eq:acq_ei}
    \alpha_{\text{EI}}(\xpar, y^*) = \mathbb{E}[(y^* - \mu(\xpar))_+] = (y^* - \mu(\xpar)) \Phi \left( \frac{y^* - \mu(\xpar)}{\sigma(\xpar)} \right) + \sigma(\xpar) \phi \left( \frac{y^* - \mu(\xpar)}{\sigma(\xpar)} \right)
    ,
\end{equation}
where $y^*$ denotes the best known value of the objective function, $\phi$ and $\Phi$ are the probability density function and the cumulative distribution function, respectively, for the standard normal.

\begin{figure}[hbt!]
	\centering
	\includegraphics[trim=0 0 0 0, clip, width=.5\textwidth]{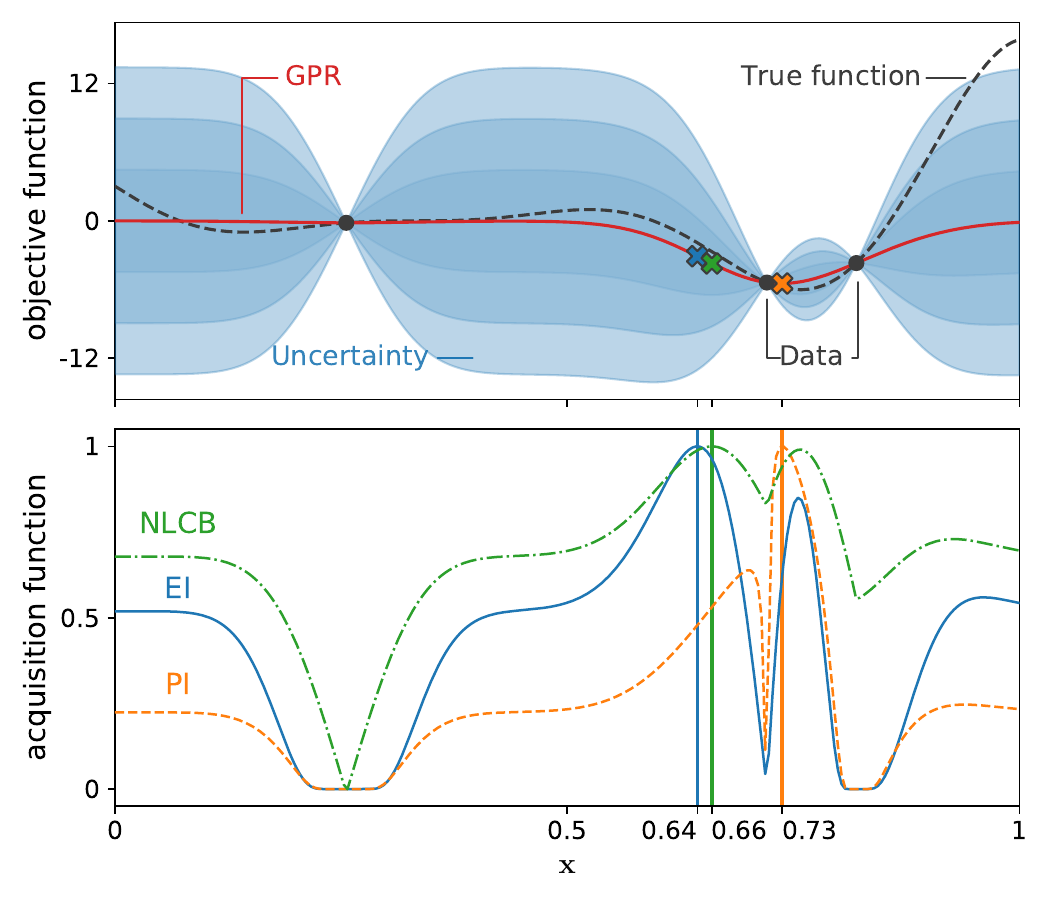}\hfill
	\includegraphics[trim=0 0 0 0, clip, width=.5\textwidth]{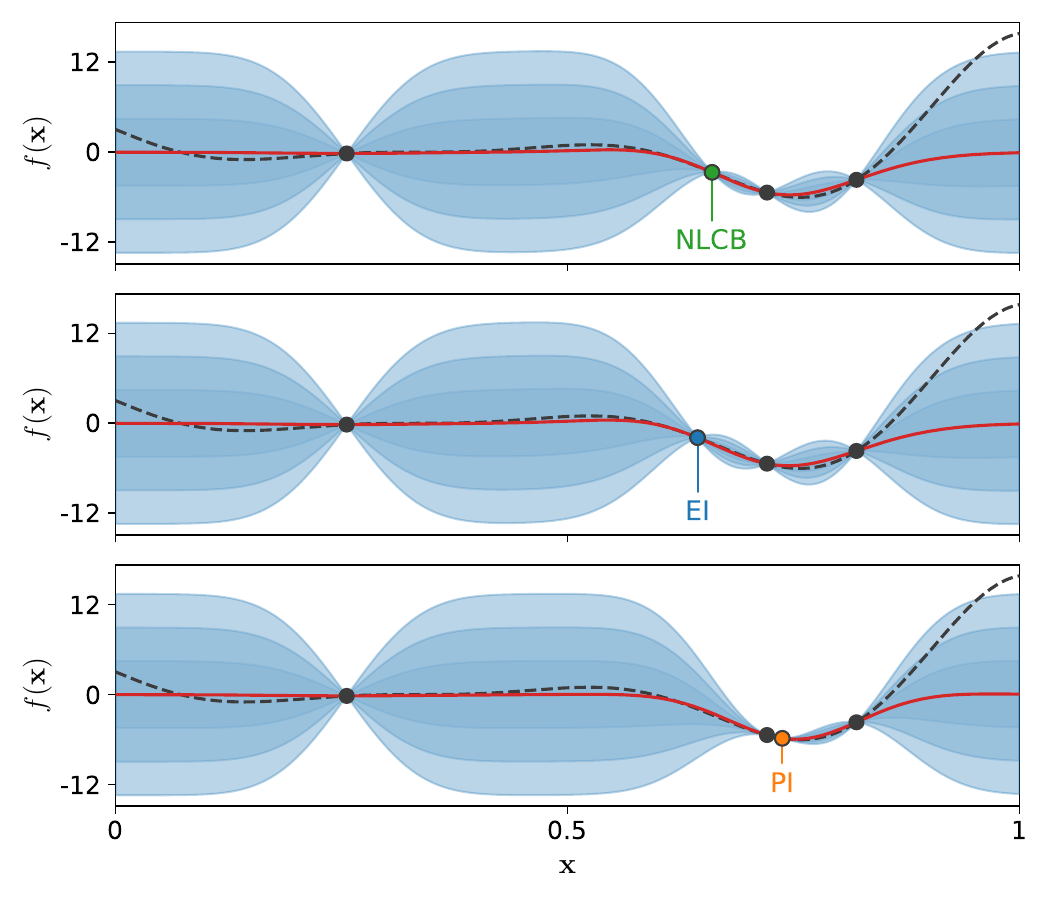}
    \caption{On the left, the GPR prediction with its uncertainty on the top and the acquisition functions at the bottom. On the right, the updated GPR after the addition of the previously selected sample, for each acquisition function.
	\label{fig:bayesian_optimization}}
\end{figure}

Probability of improvement (PI) \cite{bo_pi_1998} is an exploitative acquisition function defined as
\begin{equation}
    \label{eq:acq_pi}
    \alpha_{\text{PI}}(\xpar, y^*) = \Phi \left( \frac{y^* - \epsilon - \mu(\xpar)}{\sigma(\xpar)} \right)
    ,
\end{equation}
where $\epsilon > 0$ prevents the acquisition function from selecting candidates excessively close to the best-known solution.
A representation of a search step in BO, using different acquisition functions, is given in Figure~\ref{fig:bayesian_optimization}.

In this paper, the objective function to minimize is the physical mass of the ship.
We consider three contributions: the mass of the elements controlled by a decision variable, the mass of the elements that are fixed, and the mass of the reinforcement bars applied to correct buckling phenomena.
We add a penalization term to the physical mass to handle constraints on the QoIs that do not have a practical expression in terms of the parameters.
In particular, the number of yielded and buckled elements are required to not overly exceed some thresholds $y_{\text{crit}}$ and $b_{\text{crit}}$, respectively. Such thresholds are specified by the designers to balance the ease of optimization and the following manual post-processing of the proposed design.
The penalized objective function is obtained by summation of Equation~\eqref{eq:mass} and a penalty term:
\begin{equation}
    \label{eq:objective_mass}
    f(\xpar) = f_\text{mass}(\xpar) + f_\text{pen}(\xpar)
    .
\end{equation}
The penalization term $f_\text{pen}(\xpar)$ is defined as 
\begin{equation}
    \label{eq:penalty_mass}
    f_\text{pen}(\xpar) = c_\text{y} (n_\text{y}(\xpar) - y_{\text{crit}})_+^2 + c_\text{b} (n_\text{b}(\xpar) - b_{\text{crit}})_+^2
    ,
\end{equation}
where $c_\text{b}$ and $c_\text{y}$ are some positive constants to be chosen, chosen by the designers to quantify the trade-off between constraints violation and mass gain.
The quadratic penalization, although known to introduce non-stationarity in the BO \cite{gpr_gramacy_2021}, was preferred to more sophisticated techniques \cite{bo_constr_ineq_2014, bo_constr_ineq_2014, bo_albo_2016, gpr_gramacy_2021} since the designer can explicitly tune the trade-off between mass gain and violation.
Equation~\eqref{eq:penalty_mass} can be extended to include constraints on other QoIs, such as the vertical deflection.

We use the Emukit \cite{emukit_2019, emukit_2023} implementation of BO, where the optimization of the acquisition function is carried out by the \texttt{trust-constr} method from SciPy \cite{scipy_2020}. 
The \texttt{trust-constr} method leverages the algorithm from \cite{trust_constr_1999}, which is a barrier method that uses sequential quadratic programming and trust regions to support generic inequality constraints.
We leverage this capability to efficiently restrict the search space at runtime, by producing a bound on the linear part of the penalized objective function.

Let $\xpar^*$ be the current best solution. Then, since $n_\text{b}(\xpar) \geq 0$ and $f_\text{pen}(\xpar) \geq 0$, a better optimum can only be found in the half-space where $\xpar$ satisfies
\begin{equation}
    \label{eq:objective_ub}
    f(\xpar^*) \geq f(\xpar) \geq m_{\text{fixed}} + \mathbf{d} \cdot \xpar
    ,
\end{equation}
giving a linear inequality that will be updated each time $f(\xpar^*)$ decreases. 
This bound can be lax depending on the penalization term, but in practice it is crucial to reduce the search space and achieve good optimization performances.

The constraint on the maximum VCG allows a linear inequality formulation as well, since Equation~\eqref{eq:vcg} is linear in $\xpar$.
Then bounding from above with the critical value $\text{VCG}_{\text{crit}}$ gives
\begin{equation}
    \label{eq:vcg_ub_0}
    \text{VCG}_{\text{crit}} \geq 
    \frac{\text{VCG}_{\text{fixed}} m_{\text{fixed}} + \sum_{i} \text{VCG}_{i} \mathbf{d}_i \xpar_i}{m_{\text{fixed}} + \sum_{i} \mathbf{d}_i \xpar_i} 
    ,
\end{equation}
\begin{equation}
    \label{eq:vcg_ub_1}
    \text{VCG}_{\text{crit}} \left( m_{\text{fixed}} + \sum_{i} \mathbf{d}_i \xpar_i \right) \geq
    \text{VCG}_{\text{fixed}} m_{\text{fixed}} + \sum_{i} \text{VCG}_{i} \mathbf{d}_i \xpar_i
    ,
\end{equation}
\begin{equation}
    \label{eq:vcg_ub}
    (\text{VCG}_{\text{crit}} - \text{VCG}_{\text{fixed}}) m_{\text{fixed}} \geq
    \sum_{i} (\text{VCG}_{i} - \text{VCG}_{\text{crit}}) \mathbf{d}_i \xpar_i
    .
\end{equation}
Thus, our optimization problem is
\begin{alignat}{3}
    \label{eq:mass_opt}
    & \text{minimize} & \qquad &
        f(\xpar) = m_{\text{fixed}} + \mathbf{d} \cdot \xpar + m_{\text{bar}} n_\text{b}(\xpar) + c_\text{y} (n_\text{y}(\xpar) - y_{\text{crit}})_+^2 + c_\text{b} (n_\text{b}(\xpar) - b_{\text{crit}})_+^2
        , & \\
    & \text{subject to} & &
        f(\xpar^*) \geq m_{\text{fixed}} + \mathbf{d} \cdot \xpar
        , & \nonumber \\
    & & &
        (\text{VCG}_{\text{crit}} - \text{VCG}_{\text{fixed}}) m_{\text{fixed}} \geq \sum_{i} (\text{VCG}_{i} - \text{VCG}_{\text{crit}}) \mathbf{d}_i \xpar_i
        , & \nonumber \\
    & & & \xpar_i \in \Dpar_i & \forall i \in I 
        . \nonumber 
\end{alignat}

As an additional heuristic, we balance exploration and exploitation by switching between different acquisition functions after a number of iterations with no decrease in the current optimum.

The BO formulation used so far considers a continuous domain, but we are constrained to use commercially available thicknesses, resulting in a discrete domain.
A naive implementation would optimize the acquisition and round to the closest domain point, but this approach is prone to repeated selection of the same configuration, and does not take into account the feasibility constraint on VCG.
Multiple approaches can be found in the literature to handle ordinal and categorical variables with GPR, such as handling discretization inside the kernel~\cite{bo_discrete_trans_2020, bo_combo_2019}, while other approaches to BO replace the GPR with probabilistic models that handle discreteness transparently~\cite{bo_tpe_2011, bo_smac_2011, bo_bocs_2018}. A method capable of handling mixed variables and both known and unknown constraints can be found in~\cite{bo_mivabo_2020}.
In our application, although the domain is fully discrete, the underlying physics is inherently continuous and the QoIs exhibit discontinuities that are very small in relative value when the configuration is feasible and no penalization is being applied. 
This observation allows us to rely on rounding, after gradient-based optimization of the acquisition under known constraints, instead of using derivative-free or sampling-based acquisitions.
On detection of a duplicate sample, or if the sample is unfeasible, we solve an integer linear program (ILP) \cite{ilp_theory_1998} to find a new configuration that satisfies all constraints and has minimal distance from the original one.

Let $\bar\xpar$ be the current duplicate or infeasible candidate, $I$ be the set of parameter indices, $T$ be the set of possible parameter values, and $x_{it}$ be the binary decision variable which takes value $1$ when $\xpar_i = t$, and $0$ otherwise. 
The next sample is retrieved by solving the following optimization problem
\begin{alignat}{3}
    \label{alg:mip_rounding}
    & \text{minimize} & \qquad &
        \sum_{i \in I} \sum_{t \in \Dpar_i} x_{it} |t - \bar\xpar_i|^2 
        , & \\
    & \text{subject to} & &
        \sum_{t \in \Dpar_i} x_{it} = 1
        ,
        & \forall i \in I 
        , \nonumber \\
    & & &
        \sum_{i, t \mid \bar\xpar_i = t} x_{it} \leq | I | - 1
        , & \nonumber \\
    & & &
        \sum_{i \in I} \sum_{t \in \Dpar_i} 
         x_{it} \mathbf{d}_i t \leq m_{\text{UB}} - m_{\text{fixed}}
        , & \nonumber \\
    & & & 
        \sum_{i \in I} \sum_{t \in \Dpar_i} 
        x_{it} (\text{VCG}_{i} - \text{VCG}_{\text{crit}}) \mathbf{d}_i t \leq
        (\text{VCG}_{\text{crit}} - \text{VCG}_{\text{fixed}}) m_{\text{fixed}} 
        , & \nonumber \\
    & & & x_{it} \in \{0, 1\}, & \forall i \in I, \,\, \forall t \in T 
        , \nonumber 
\end{alignat}
where $\mathbf{d}_i$ is the linear density of the $i$-th parameter and $m_{\text{UB}}$ the penalized mass of the current optimum.
The first constraint ensures consistency, the second one guarantees that the new sample will not be a duplicate of $\bar{\xpar}$, and the rest replicate the constraint of the BO.
If the newly found candidate is still a duplicate, we apply a random disturbance until an actual unvisited configuration is generated, or a maximum number of attempts is reached. 
For the solution of the ILP, we use the COIN-OR Branch-and-Cut solver (CBC) \cite{cbc_2023} through the python-mip library~\cite{python_mip_2020}. 
This ILP has a small size, comprising only $|T||I|$ integer variables and $|I|+3$ constraints, and is often solved to optimality by the preprocessing of the solver.

BO selects new candidate configurations until a computational budget is expended, given as the number of iterations or total execution time. 
The complete procedure is reported in Algorithm~\ref{alg:bo}.

\begin{algorithm}[ht]
\begin{algorithmic}
    \Require objective function $f$, high-fidelity parameters $\Xpar$, and evaluations $\ypar$
    \Ensure returns the configuration that minimizes $f$
    \State train GPR on $\Xpar$, $\ypar$
    \State $\xpar^* \gets \argmin_{\xpar \in \Xpar} f(\xpar)$
    \State $\alpha \gets \text{one of } \{\alpha_{\text{NLCB}}, \alpha_{\text{EI}}, \alpha_{\text{PI}}\}$
    \While{number of iterations $\leq$ iter limit}
        \State $\xpar \gets \argmax_{\xpar} \alpha(\xpar)$ subject to Equation~\eqref{eq:objective_ub}, Equation~\eqref{eq:vcg_ub}
        \State $\xpar \gets \text{round}(\xpar)$ \Comment{could be provided by the BO framework}
        \While{$\xpar$ is infeasible \textbf{or} $\xpar \in \Xpar$}
            \State $\xpar \gets $ solution to ILP~\eqref{alg:mip_rounding}
        \EndWhile
        \If{$f(\xpar) < f(\xpar^*)$}
            \State $\xpar^* \gets \xpar$
            \State update $f(\xpar^*)$ in Equation~\eqref{eq:objective_ub}
        \EndIf
        \State $\Xpar \gets \Xpar \cup \{\xpar^*\}$
        \State $\ypar \gets \ypar \cup \{f(\xpar^*)\}$
        \State train GPR on $\Xpar$, $\ypar$
        \If{$\ypar^*$ did not decrease in the last 100 iterations}
            \State $\alpha \gets \text{another one of } \{\alpha_{\text{NLCB}}, \alpha_{\text{EI}}, \alpha_{\text{PI}}\}$
        \EndIf
    \EndWhile \\
    \Return $\xpar^*$ 
\end{algorithmic}
    \caption{Specialized BO implementation.}
    \label{alg:bo}
\end{algorithm}

If the search finds no candidate that improves on the current best solution, the optimization is concluded. 
Otherwise, a high-fidelity simulation is required to validate such result and the search is repeated. 
In practice, even with the available refinements to reduce the search space, BO on a large parameter space is only able to perform several rounds before failing to find new promising candidates.

\subsubsection{Principal dimensions search}
\label{subsubsec:opt_pds}
To overcome the limitations of repeated BO on our large parameter space, we propose a more principled greedy heuristic, based on the exhaustive exploration of the neighborhood of the best candidate configuration.
By leveraging the cheap surrogates, it is possible to evaluate all configurations that differ for one parameter from the chosen candidate.
This approach is related to the cyclic coordinate search but performs a full scan of each parameter, thus the name principal dimensions search (PDS). 
The maximum number of surrogate evaluations is $\mathcal{O}(\sum_{i \in I}|\mathcal{D}_i|-1)$, but in practice it is much lower due to the enforcement of Equation~\eqref{eq:objective_ub} and Equation~\eqref{eq:vcg_ub}. 
If no improvement is found, the search is concluded. 
Otherwise, the best configuration is used as the starting point of a new PDS until a computational budget is exhausted, and the best results are validated by the high-fidelity solver.
The entire process is detailed in Algorithm~\ref{alg:pds}.

\begin{algorithm}[ht]
\begin{algorithmic}
    \Require objective function $f$, high-fidelity parameters $\Xpar$
    \Ensure returns a feasible solution through PDS
    \State $\xpar^* \gets \argmin_{\xpar \in \Xpar} f(\xpar)$
    \While{number of iterations $\leq$ iter limit \textbf{and} elapsed time $\leq$ time limit}
        \State $\xpar^{\text{base}} \gets \xpar^*$
        \State $\xpar \gets \xpar^{\text{base}}$
        \For{i in I}
            \For{$t \in \Dpar_i  \setminus \{ \xpar^{\text{base}}_i \}$}
                \State $\xpar_i \gets t$
                \If{$\xpar$ is feasible \textbf{and} $f(\xpar) \leq f(\xpar^*)$} 
                    \State $\xpar^* \gets \xpar$
                \EndIf
            \EndFor
            \State $\xpar_i \gets \xpar^{\text{base}}_i$
        \EndFor
    \EndWhile\\
    \Return $\xpar^*$, $f(\xpar^*)$
\end{algorithmic}
    \caption{Principal direction search.}
    \label{alg:pds}
\end{algorithm}
This greedy heuristic is strongly exploitative of the current best solution and proves effective in finding new candidates but is limited in its exploratory capabilities.

\subsection{Parameterization refinement}
\label{subsec:reparam}
The results of the optimization procedure are determined by the problem formulation, that is, by the parameterization chosen for the model.
While a human designer is able to, in principle, assign the thickness of each element independently, the optimization of a model with so many decision variables would be impractical even for problems of moderate scale.
On the other hand, a model that only uses a few parameters will be limited in its ability to balance the optimization of total mass and the reduction of structural failure phenomena, which are often quite local to critical geometrical features of the mesh.
A first simplification comes from observing that, in the shipyard, the structure of the ship is assembled by metal sheets with dimensions about \SI{2.5}{m}$\times$\SI{15}{m}, which correspond to a group of about 90 elements, on average. 
The elements that are part of the same sheet form a patch, and this is the finest level of parameterization that the designers can use in formulating the optimization problem.
On the other hand, the coarsest parameterization consists of assigning a decision variable for each group of patches that share the same regulatory minimum thickness.
The common practice is to also take into consideration the orientation of the patches and the possible criticalities emerging from the stress configurations under the load conditions.
Upon convergence of the optimization pipeline, the designer analyzes the optimum configuration to identify whether it is possible to improve upon its results.
A common situation is that a large group of patches, controlled by the same decision variable, is assigned a high thickness even if only a small, localized subgroup of those patches would incur in yielding and buckling.
The designers would then split the original grouping of the patches, to manually assign different thickness values to the subgroups.
A lower thickness would be assigned to the subgroup that does not incur in structural issues, and higher thickness to the one where the failure phenomena are localized.
This procedure could result in lowering the total mass of the model, depending on the difference in thickness and the extent of the structural criticalities.
It would be natural to derive a new optimization problem where the two subgroups of patches are assigned to two different decision variables, but this approach would incur two expensive issues: the setup of a new problem is a cumbersome operation, and the new optimization would require a new initial sampling and surrogate construction.

We propose an automatic refinement of the parameterization based on the solution of a set of ILPs, which will determine a division of each group of elements according to its structural response.
This approach has been inspired by a similar idea in the context of multiscale finite elements \cite{clustering_benaimeche_2022}, where the number of finer, expensive problems to be solved was reduced by clustering the boundary conditions applied to the sub-elements.

In our case, we start from the best solution obtained during the optimization and the set of surrogates constructed based on the high-fidelity evaluations collected during the pipeline execution.
A number of configurations are evaluated in the vicinity of the best solution, to obtain the yielding and buckling states of each element.
For each parameterized section, an ILP is constructed to choose the best parameter value for each patch, so that an objective function similar to the one used in BO is minimized.
The resulting thickness assignments identify the optimal refined parameterization.

Let $P_i$ be the set of patch indices assigned to the $i$-th parameter $\xpar_i$, and $\Dpar_i$ the set of possible parameter values. 
Additionally, let $y_{pt}$ and $b_{pt}$ be the number of yielded and buckled elements, respectively, in patch $p \in P_i$ when $\xpar_i = t$. 
We assume that each patch is independent of the rest of the ship, so that $y_{pt}$ and $b_{pt}$ combine linearly and can be computed (using surrogates, and high-fidelity evaluations when available) by enumerating $t \in T$ and keeping the rest of the parameter values fixed.
These evaluations follow the same order as in Algorithm~\ref{alg:pds}.

The binary decision variables $x_{pt}$ express the assignments of patches to parameter values, with $x_{pt} = 1$ only if patch $p$ has thickness $t$. 
For consistency, for each patch $p$, it is required that exactly one of the $x_{pt}$ variables is nonzero.
The upper bound on the VCG takes the same form used in Equation~\eqref{eq:vcg_ub}. 
For the objective function, we aim to minimize the total mass of the patches in $P_i$, their reinforcement bars, and penalty terms given by the sum of squares of $y_{pt}$ and $b_{pt}$. 
Finally, we fix the actual number of different parameter values $n_{\text{clusters}} \geq 2$ through the additional decision variables $u_t$ and its consistency constraints.

We obtain the following optimization problem:
\begin{alignat}{3}
    \label{alg:mip_clustering}
    & \text{minimize} & \qquad &
        \sum_{p \in P_i} \sum_{t \in \Dpar_i} 
        x_{pt} (\mathbf{d}_p t + m_\text{bar} b_{pt} + c_\text{y} y_{pt}^2 + c_\text{b} b_{pt}^2 )
        ,
        & \\
    & \text{subject to} & &
        \sum\limits_{t \in \Dpar_i} x_{pt} = 1, & \forall p \in P_i 
        , \nonumber \\
    & & & x_{pt} \leq u_{t}, & \forall p \in P_i, \,\, \forall t \in \Dpar_i 
        , \nonumber \\
    & & & \sum_{p \in P_i} x_{pt} \geq u_{t}, & \forall t \in \Dpar_i 
        , \nonumber \\
    & & & \sum_{t \in \Dpar_i} u_{t} = n_{\text{clusters}} , & \nonumber \\
    & & & \sum_{p \in P_i} \sum_{t \in \Dpar_i} 
        x_{pt} (\text{VCG}_{p} - \text{VCG}_{\text{crit}}) \mathbf{d}_p t
        \leq (\text{VCG}_{\text{crit}} - \text{VCG}_{\text{res}}) m_{\text{res}} , & \nonumber \\
    & & & x_{pt} \in \{0, 1\}, & \forall p \in P_i, \,\, \forall t \in \Dpar_i 
        , \nonumber \\
    & & & u_{t} \in \{0, 1\}, & \forall t \in \Dpar_i
        , \nonumber 
\end{alignat}
where $\mathbf{d}_p$ is the linear density of patch $p$, $\text{VCG}_p$ is the VCG of patch $p$, and $\text{VCG}_{\text{res}}$ and $m_{\text{res}}$ are the VCG and mass, respectively, of all elements not controlled by the current parameter.

The total number of integer variables is $(|P_i|+1)|\Dpar_i|$ and the number of constraints is $|P_i||\Dpar_i|+|P_i|+|\Dpar_i|+2$.
This class of problems is NP-hard, but branch and cut simplex-based solvers are extremely efficient in practice thanks to pre-solve transformations and branching rules~\cite{mip_simplex_2004, mip_12years_2013}.
An example of the procedure, for a simplified section composed of 6 patches, is depicted in Figure~\ref{fig:reparam_simple_06_flow}.

\begin{figure}[hbt!]
	\centering
	\includegraphics[trim=0 0 0 0, clip, width=1\textwidth]{figures/reparam_simple_06.pdf}
    \caption{The reparameterization procedure on a simplified parameterized section. The section is composed of 6 patches, represented in steps 1, 2, and 3 with different colors corresponding to different thicknesses. In step 1, all patches have the same thickness as selected by the previous optimization. In step 2, different thicknesses are evaluated to obtain the structural responses of each patch. In step 3, the structural responses determine the reparameterization problem, and its optimal solution clusters the patches in two groups. In step 4, the patches are assigned the optimized thicknesses, and two parameterized sections are determined.
	\label{fig:reparam_simple_06_flow}
    }
\end{figure}

We remark that the objective functions in ILP~\eqref{alg:mip_clustering} and Equation~\eqref{eq:objective_mass} are fundamentally different in the way they handle the penalization of yielded and buckled elements, as the former uses a sum of squares and the latter uses the square of a sum.

In ILP~\eqref{alg:mip_clustering}, the objective function can be split in the contributions from each patch to the physical mass, due to thickness and additional stiffeners, and the penalization of failed elements as
\begin{equation}
    \sum_{p \in P_i} \sum_{t \in \Dpar_i} 
        x_{pt} (\underbrace{\mathbf{d}_p t \vphantom{y_{pt}^2}}_{\text{patches mass}} + \underbrace{m_\text{bar} b_{pt} \vphantom{y_{pt}^2}}_{\text{reinforcement bars}} + \underbrace{c_\text{y} y_{pt}^2 + c_\text{b} b_{pt}^2}_{\text{per-patch penalty}} )
    ,
    \label{eq:mip_linear_obj}
\end{equation}
which is linear in the decision variables $x_{pt}$.
The proper reformulation of Equation~\eqref{eq:objective_mass} in terms of $x_{pt}$ requires instead the quadratic programming expression
\begin{equation}
    \sum_{p \in P_i} \sum_{t \in \Dpar_i} 
        x_{pt} (
        \underbrace{\mathbf{d}_p t\vphantom{y_{pt}^2}}_{\text{patches mass}}
        + \underbrace{m_\text{bar} b_{pt}\vphantom{y_{pt}^2}}_{\text{reinforcement bars}}
        )
    + f_\text{pen}(x)
    ,
    \label{eq:mip_quadratic_obj}
\end{equation}
where the penalty term for the binary decision variables, corresponding to Equation~\eqref{eq:penalty_mass}, is
\begin{equation}
    f_\text{pen}(x) = \underbrace{
        \left( \left( \sum_{p \in P_i} \sum_{t \in \Dpar_i} x_{pt} c_\text{y} y_{pt} \right) - y_\text{crit} \right)_+^2}_{\text{global yielding penalty}}
        + \underbrace{
        \left( \left( \sum_{p \in P_i} \sum_{t \in \Dpar_i} x_{pt} c_\text{b} b_{pt} \right) - b_\text{crit} \right)_+^2}_{\text{global buckling penalty}}
    .
    \label{eq:mip_quadratic_penalty}
\end{equation}

The thresholds $y_\text{crit}$ and $b_\text{crit}$ which in Equation~\eqref{eq:mip_quadratic_penalty} need to be adapted to the current $i$-th parameterized section, as Equation~\eqref{eq:penalty_mass} uses global values.
One option would be subtracting the yielded and buckled elements totals of the other sections from the critical values, but the combination of the separate problems could still exceed the global constraints.
As an alternative, ILP~\eqref{alg:mip_clustering} needs to be reformulated as a global optimization problem, optimizing all the parameterized sections at once, which is impractical for large models.
In practice, the usage of ILP~\eqref{alg:mip_clustering} leads to a more severe penalization of failure phenomena in each patch, so that the updated model will offer better trade-offs during the optimization.

By solving ILP~\eqref{alg:mip_clustering} for each parameterized section, and possibly multiple $n_{\text{clusters}}$, a collection of candidate reparameterizations is obtained.
At most one refinement can be applied to each parameterized section, including the trivial one which leaves it unchanged.
Moreover, the designers could specify a maximum number of parameters for the model.
The selection of the optimal set of refinements can be obtained by the solution of a knapsack ILP \cite{ilp_theory_1998}, where the cost of a refinement is the number of clusters and its value is given by ILP~\eqref{alg:mip_clustering}.

Since the construction of the VCG constraint assumes that all the elements outside the current group of patches do not change mass, the parameter configuration obtained by joining the optimized assignments might violate the global VCG constraint.
Although an infeasible super-optimal configuration is undesirable in principle, the goal of the procedure is to obtain an expressive parameterization, that will be leveraged by the following optimization.

Once the optimal selection of clusters has been chosen, the $i$-th parameterized section associated to $\xpar_i$ is also partitioned in a set of patch clusters $\left\{ \mathcal{P}_j \right\}_{j=1}^{m}$ with $1 \leq m \leq n_\text{clusters}$.
The patches in each cluster $\mathcal{P}_{j > 1}$ will be assigned to a new parameter $\xpar_{k}$, with feasible domain $\mathcal{D}_i$, decreasing the number of patches assigned to $\xpar_i$.
With the new parameterization, the configurations generated with the coarse model can still be retrieved by enforcing the constraints 
\begin{equation}
    \xpar_i = \xpar_{k}\qquad\forall k \in \mathcal{K}_i
    ,
    \label{eq:reparam_constraint}
\end{equation}
where $\mathcal{K}_i$ collects the indices of the $m-1$ parameters generated for the clusters $\mathcal{P}_{j > 1}$.
This feature is crucial for the efficiency of the optimization procedure, as all the previous high-fidelity simulations need not be executed again.
On the other hand, all the high-fidelity samples collected so far can be viewed as coming from a constrained version of the current parameter space, where Equation~\eqref{eq:reparam_constraint} reduced the feasible configurations.
Due to this consideration, the surrogates must be rebuilt on a sampling that includes a configuration where Equation~\eqref{eq:reparam_constraint} is not active, so that the effects of the new parameterization are properly observed.
For this task, we repeatedly perform a random sampling of the new parameter space, keeping the candidate set which maximizes the cumulative distance between the high-fidelity samples and the new ones.

At this point, the surrogates are rebuilt on the updated high-fidelity samples and the optimization, both multi- and single-objective, can take place on a more expressive model.
With this approach, the parameterization is adapted as the optimization progresses and is thus able to overcome the biases that the designers could have introduced in the initial model creation.
This presents a critical advantage when the hull features high complexity due to novel structures or design constraints, and enables the pipeline to provide high-quality initial designs with minimal supervision by the users.

\section{Numerical results}
\label{sec:results}
This section presents the results from the application of our optimization pipeline to two ship models: a simplified midship section typical of the initial design phase in Section~\ref{subsec:results_msection}, and a full ship model with all environments and functional features in Section~\ref{subsec:results_fullship}.

\subsection{Midship section}
\label{subsec:results_msection}
\begin{figure}[hbt!]
	\centering
	\includegraphics[trim=0 0 0 0, clip, width=.33\textwidth]{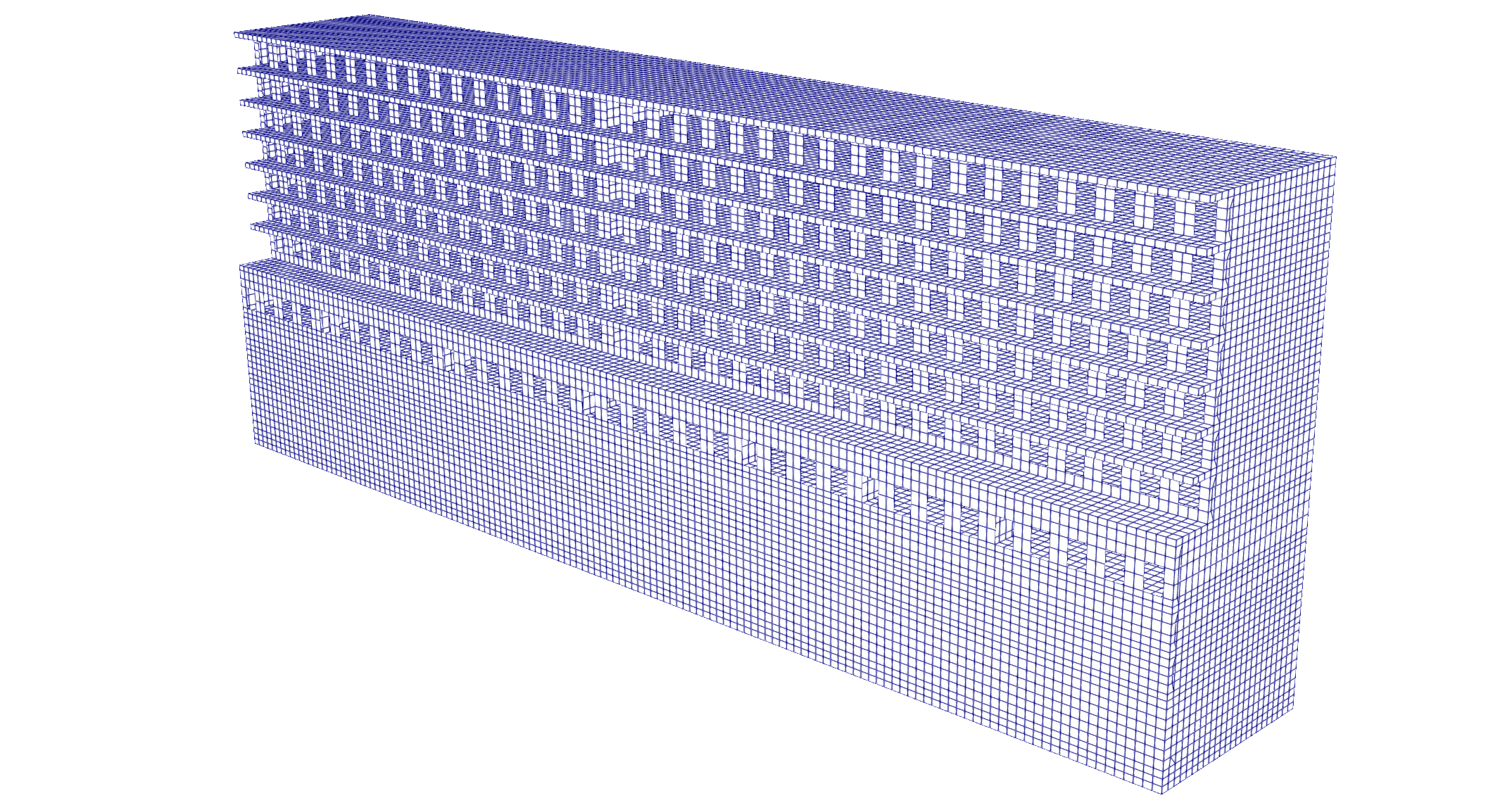}\hfill
	\includegraphics[trim=0 0 0 0, clip, width=.33\textwidth]{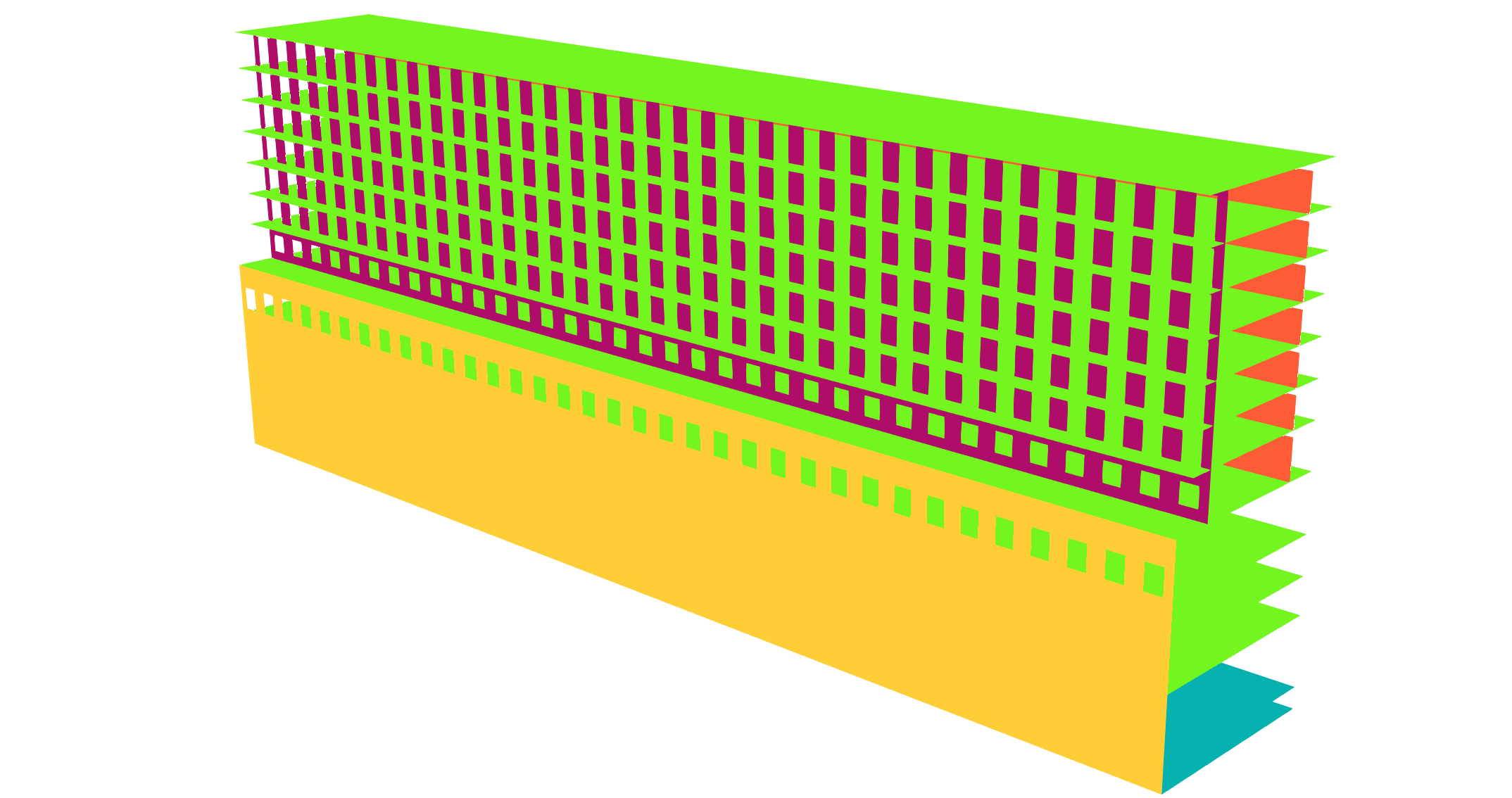}\hfill
	\includegraphics[trim=0 0 0 0, clip, width=.33\textwidth]{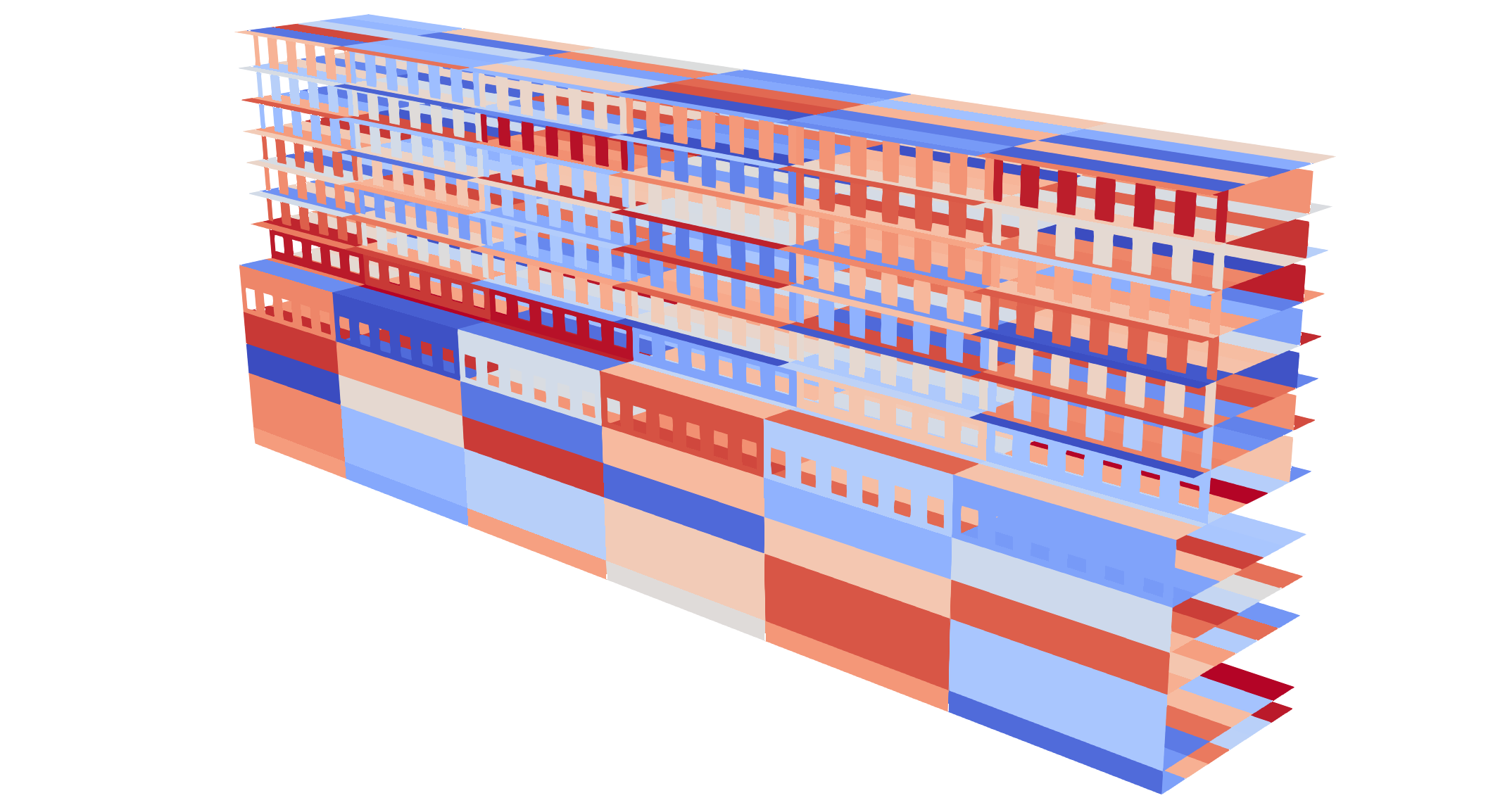}
	\caption{Full view of the midship section on the left, the initial parameterization in the middle, and the details of the parameterized patches on the right.}
	\label{fig:msection_general}
\end{figure}

For the initial development and testing of our methods, we choose a simplified model of the main section of a typical hull. 
The model is depicted in Figure~\ref{fig:msection_general}. 
It consists of $75192$ elements, of which, $52360$ are assigned to a parameterized section as reported in Table~\ref{tab:msection_params}. 
The parameterized elements are further grouped in $582$ patches.
The initial parameterization follows the orientation of the patches and the regulatory minimum thickness of the structural members.
Thus, lower bounds for the parameter values come from standard regulations, while the upper bounds have been chosen according to common design practices. The total number of configurations in the domain is $62720$.
This model is a quarter of the actual hull: during the high-fidelity simulations, it is reflected across the $xz$ and $yz$ planes. 

A NASTRAN run on this model takes about $1$ minute and the query time of the surrogates is $0.15$ seconds, for a speedup factor of $400$. 
The maximum value for the VCG is \SI{15}{m}. 
In this experiment, the critical values for the yielding due to the stresses and the von Mises criterion are lower than the regulatory standard, as the simplicity of this model would otherwise result in no yielded elements. 
For this model, the maximum number of parameters is $20$.

\begin{table}[hbt!]
	\caption{Parameters description of the midship section test case. Thickness values are in \si{mm}.}
    \label{tab:msection_params}
	\centering
	\begin{tabular}{ c l r r r r r }
		\hline
		\hline
		Parameter & Region & Patches & Default & Min. & Max. & \# values \\
		\hline
		\hline
		\rowcolor{Gray}
		$\xpar_{1}$ & Bottom and inner bottom    & 72  & 14.0 & 12.0 & 20.0 & 14 \\
		$\xpar_{2}$ & Decks from $2$ to $12$     & 396 & 5.0  & 5.0  & 15.0 & 5 \\
		\rowcolor{Gray}
		$\xpar_{3}$ & External bulkheads         & 30  & 10.0 & 8.0  & 15.0 & 14 \\
		$\xpar_{4}$ & Internal bulkheads         & 42  & 5.0  & 5.0  & 15.0 & 8 \\
		\rowcolor{Gray}
		$\xpar_{5}$ & Shell plating              & 42  & 8.0  & 8.0  & 15.0 & 8 \\
		\hline
		\hline
	\end{tabular}
\end{table}

We start the optimization pipeline by performing a random sampling of the parameter space, selecting $20$ configurations in addition to the default one. 
We remark that the choice of a small initial sampling is motivated by the subsequent multi-objective optimization, in which high-fidelity experiments are performed on parameter configurations sampled by the surrogate Pareto frontier.
The truncation rank chosen for the POD is 6, so that the discarded modes have normalized module smaller than $0.01$.
This is a relatively high ratio compared with the POD literature, but finds its justification in that we are not interested in the stress field per se, but in sufficiently accurate modeling of the number of yielded and buckled elements.
In Appendix~\ref{subsec:appendix_midship}, the truncation rank and the resulting errors on the QoIs are discussed in detail.
The GPRs for the reduced coefficients employ the squared exponential kernel with ARD and epistemic noise variance for numerical stability.

\begin{figure}[hbt!]
	\centering
	\includegraphics[trim=0 0 0 0, clip, width=\textwidth]{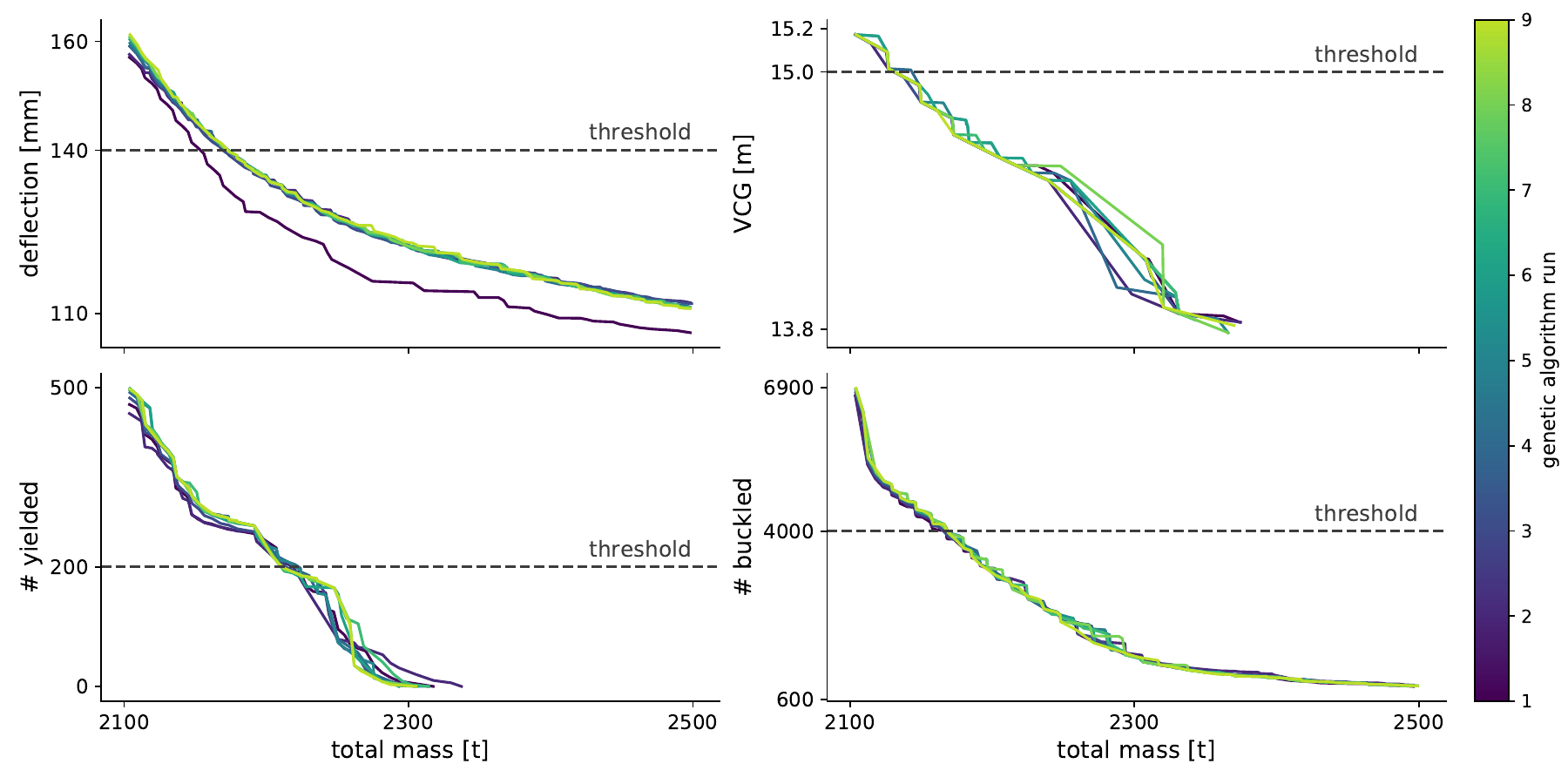}
	\caption{Evolution of the surrogate Pareto frontier in subsequent runs (generations) of the genetic algorithm, for the midship section with 5 parameters. 
    Each figure depicts the 2-dimensional frontier computed from the projection of the full 5-dimensional frontier.
    }
	\label{fig:msection_pareto_evolution}
\end{figure}

The multi-objective optimization enables the designers to inspect the effective trade-offs between different QoIs, so that the critical thresholds for the constraints can be determined from a global perspective.
Figure~\ref{fig:msection_pareto_evolution} shows the progress of the surrogate PF as the pipeline repeats a sequence of genetic algorithm run, infill selection and high-fidelity validation, and surrogates update.
The population size is set to 2000 configurations, the number of generations is 10, and the number of samples selected by the infill criterion is 9.
For the midship section, the acceptable numbers of yielded and buckled elements, and the maximum vertical deflection, are set to 200, 4000, and 140~\si{mm}, respectively.

\begin{figure}[htb]
	\centering
	\includegraphics[trim=0 0 0 0, clip, width=\textwidth]{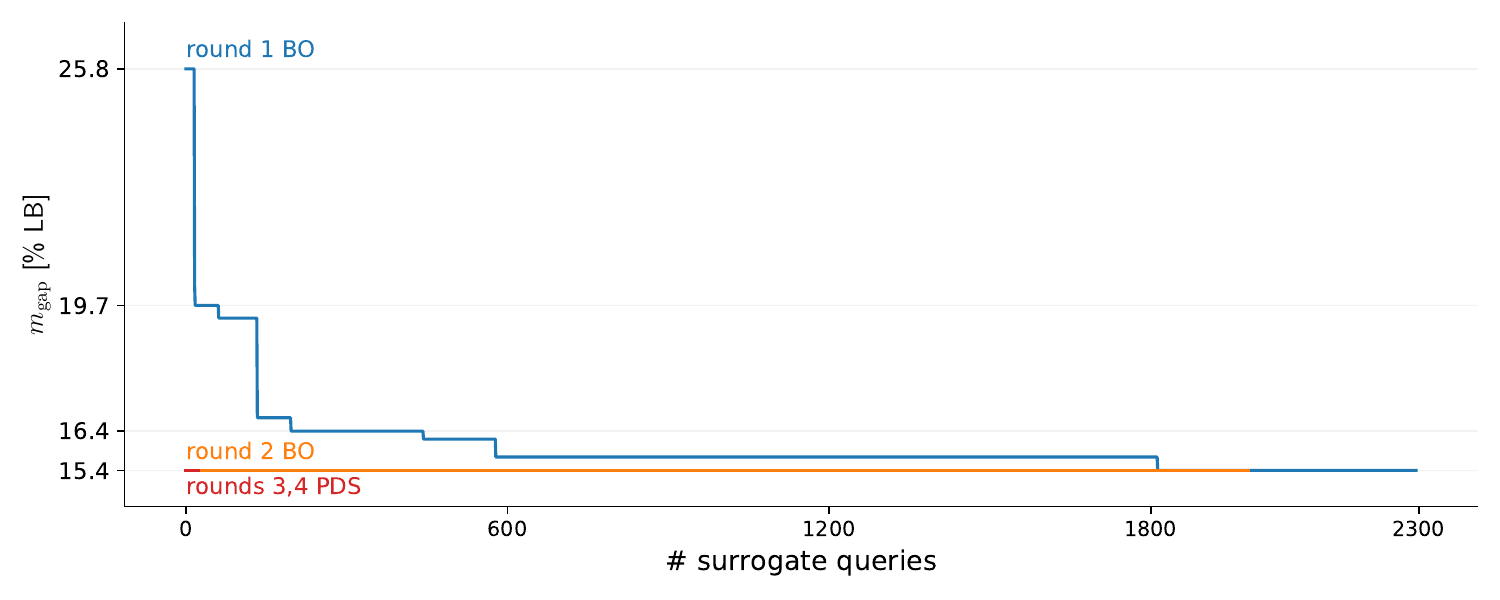}
	\caption{History of the surrogate best solution during the mass optimization of the midship section, 5 parameters. After the end of each round, the best configurations found are validated by a high-fidelity simulation, and the surrogates are rebuilt on the updated database.}
	\label{fig:msection_opt_mass_5p_history}
\end{figure}

Single-objective mass optimization starts with BO and a time limit of 5 minutes.
The history of the best surrogate candidates is presented in Figure~\ref{fig:msection_opt_mass_5p_history} in terms of the percentage gap from the theoretical lower bound (LB), that is
\begin{equation}
    m_{\text{gap}}(\xpar) = 100 \frac{\mathbf{d} \cdot (\xpar - \xpar_{\text{LB}}) + m_{\text{bar}} n_\text{b}(\xpar) + f_\text{pen}(\xpar) }{\mathbf{d} \cdot \xpar_{\text{LB}}} 
    .
    \label{eq:mass_gap}
\end{equation}
This quantity highlights the effect on the mass actually controlled by the parameters and does not depend on how extensively the model has been parameterized.
In this case, $m_{\text{fixed}} = 1108.21$~\si{t} and $\mathbf{d} \cdot \xpar_{\text{LB}} = 975.55$~\si{t}.

Only the first round is able to find better candidates, which are then confirmed by the following high-fidelity experiment. 
The second BO execution does not improve on the previous, nor do the subsequent PDS refinements. 

\begin{figure*}[hbt!]
	\centering
	\includegraphics[trim=0 0 0 0, clip, width=.5\textwidth]{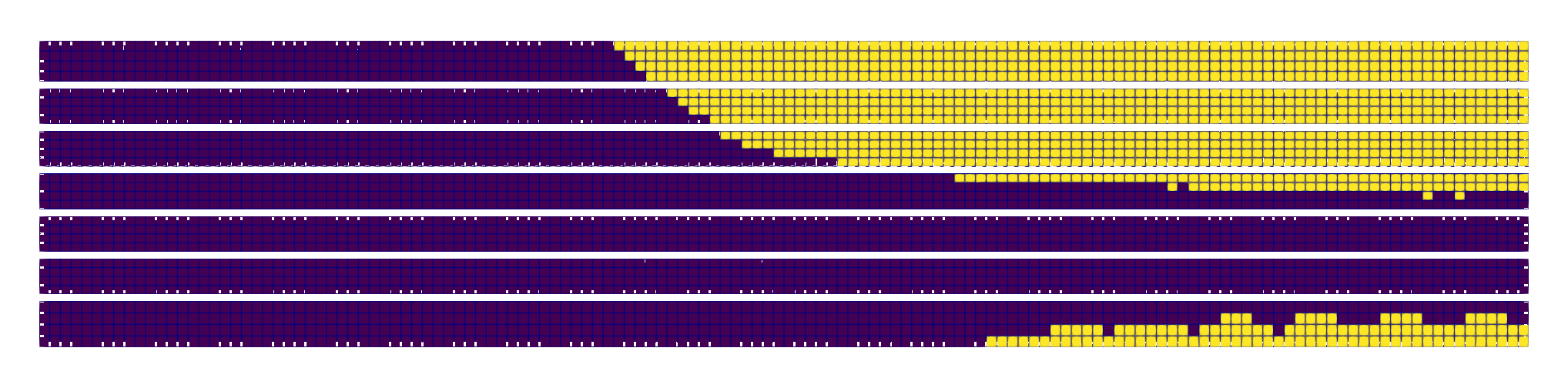}\hfill
	\includegraphics[trim=0 0 0 0, clip, width=.5\textwidth]{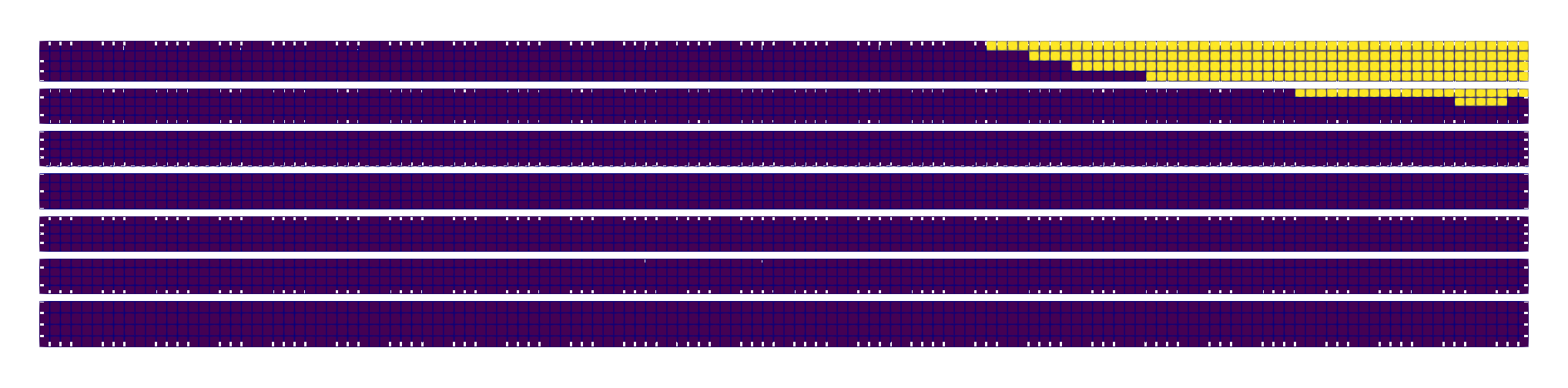}
	\caption{Buckling state of the elements controlled by $\xpar_{4}$, with failed elements in yellow and healthy elements in purple. The state on the left is the state obtained for $\xpar_{4} = 6~\si{mm}$, the right one for $\xpar_{4} = 9~\si{mm}$.}
	\label{fig:msection_p04_buckled_6_9}
\end{figure*}

\begin{figure*}[hbt!]
	\centering
	\includegraphics[trim=0 0 0 0, clip, width=.70\textwidth]{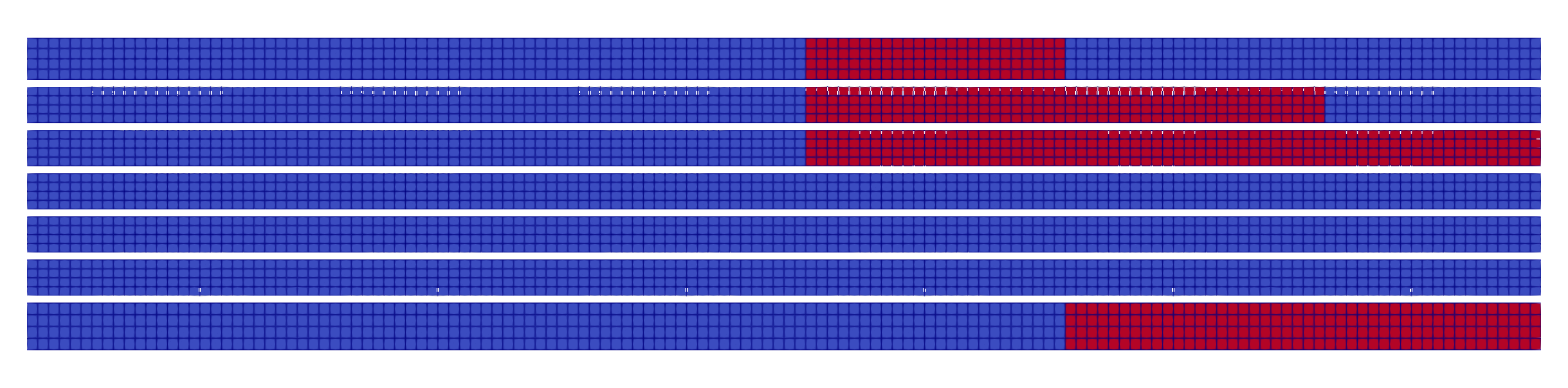}
	\caption{Optimal clustering of the elements controlled by $\xpar_{4}$. The elements in red will be associated to a new decision variable, while the elements in blue remain associated to $\xpar_{4}$}
	\label{fig:msection_p04_clusters_10p}
\end{figure*}

At this point, the optimum configuration is analyzed for the parameterization refinement, as in Section~\ref{subsec:reparam}.
For each parameter, the number of possible clusters is set to $2$.
In Figures~\ref{fig:msection_p04_buckled_6_9} and \ref{fig:msection_p04_clusters_10p}, we focus on the internal bulkheads controlled by $\xpar_4$.
The optimal clustering selects thicknesses $6$ and $9$~\si{mm}, which when applied to the entire parameter give the buckling patterns in Figure~\ref{fig:msection_p04_buckled_6_9}.
Finally, Figure~\ref{fig:msection_p04_clusters_10p} shows the assignment of elements to the two clusters.
Notably, the $3$ patches in the top right corner were assigned a lower thickness value: the ILP solution leverages the trade-off in using reinforcement bars instead of thicker, heavier patches.


After the update of the parameterization, the situation is analogous to an initial problem formulation, but with many high-fidelity configurations already available.
However, these configurations were generated from a coarser parameterization, and are not representative of the interactions between the newly created parameters.
To create the surrogate models on the new parametric domain, we generate $20$ configurations by random sampling from a uniform distribution on the updated domain.
We produce multiple candidate sets and select the one maximizing the minimum distance between the candidates, and between candidates and high-fidelity samples.
The candidate set is passed to the high-fidelity solver to generate the corresponding snapshots, and the surrogates for the updated parametric domain are created.
At this point, the pipeline resumes with multi-objective and single-objective optimization as done previously, but with a higher number of decision variables and a more expressive problem formulation.
The sequence of alternating optimization and parameterization refinement can be repeated until a set number of parameters is reached, or the prospective decrease in objective function incurs in diminishing returns.

\begin{figure}[hbt!]
	\centering
	\includegraphics[trim=0 0 0 0, clip, width=\textwidth]{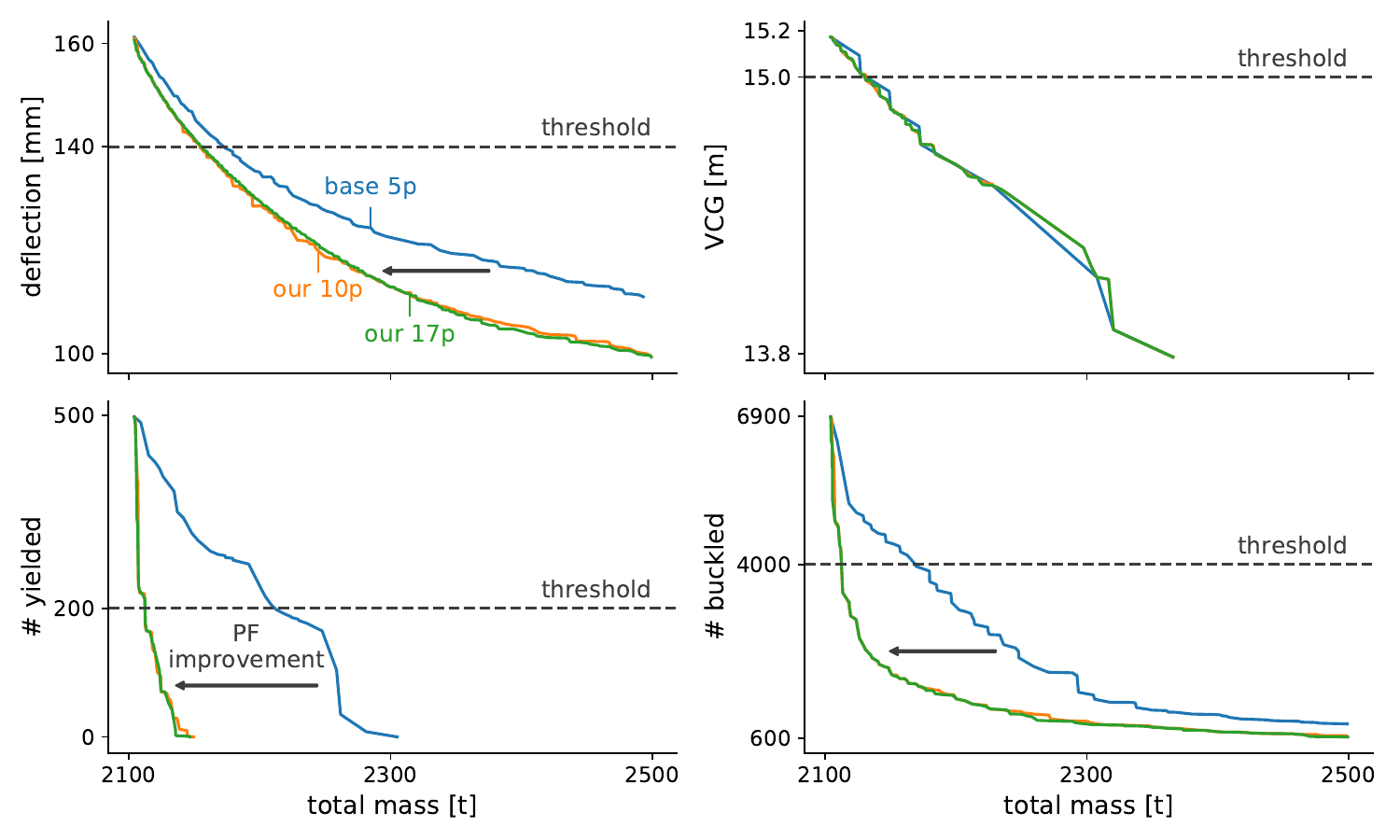}
	\caption{Evolution of the surrogate Pareto frontier for the midship section, for different parameterizations.}
	\label{fig:msection_pareto_frontiers}
\end{figure}

\begin{figure}[htb]
	\centering
	\includegraphics[trim=0 0 0 0, clip, width=\textwidth]{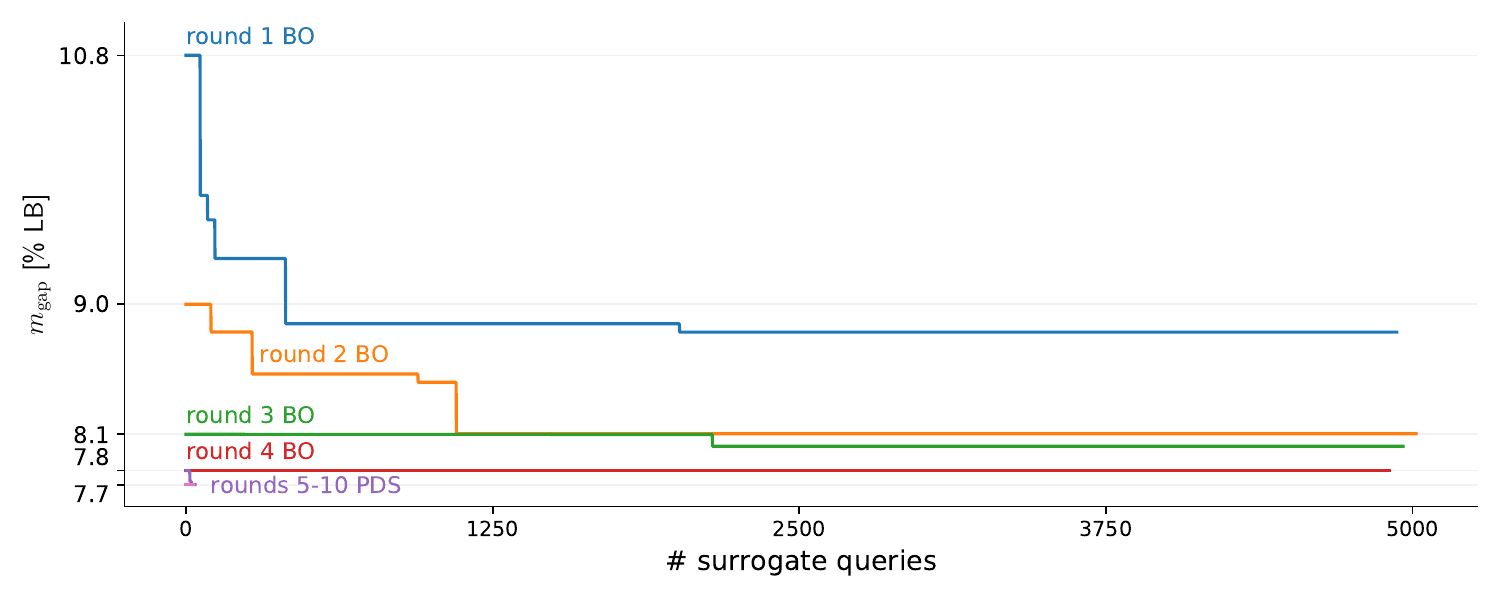}
	\caption{History of the surrogate best solution during mass optimization of the midship section with 10 parameters. After the end of each round, the best configurations found are validated by a high-fidelity simulation, and the surrogates are rebuilt on the updated database.}
	\label{fig:msection_opt_mass_10p_history}
\end{figure}

For the midship section, we perform two reparameterization steps, the first creating $5$ new parameters for a total of $10$, and the second creating $7$ new parameters for a total of $17$.
Each time, the truncation rank of the POD is increased by the number of added parameters.
Figure~\ref{fig:msection_pareto_frontiers} shows the differences in the final Pareto frontiers as the number of parameters increase.
The largest improvement in PF quality is attributed to the first reparameterization, as seen in the plots for the number of yielded and buckled elements.
The updated PF intercepts the thresholds at a much lower total mass value, indicating that the subsequent single-objective optimizations could find much lighter configurations with no additional drawbacks.
The second refinement does not produce an appreciable change in the PFs and thus the procedure is stopped.
Figure~\ref{fig:msection_opt_mass_10p_history} shows the single-objective optimization of the $10$ parameters problem.
The finer parameterization enables the optimizer to find better configurations, with the PDS being repeated several times.
In this case, the high-fidelity validation disproves a large number of optimum candidates.

We also test the effect of a single reparameterization exhausting the parameters budget, which we refer to as "one-shot", by generating $15$ new parameters from the optimal 5-parameters configuration.
Finally, we optimize a $20$ parameters model produced by the designers to compare the efficacy of our method.
The POD truncation ranks for both cases is $21$.

\begin{figure}[hbt!]
	\centering
	\includegraphics[trim=0 0 0 0, clip, width=\textwidth]{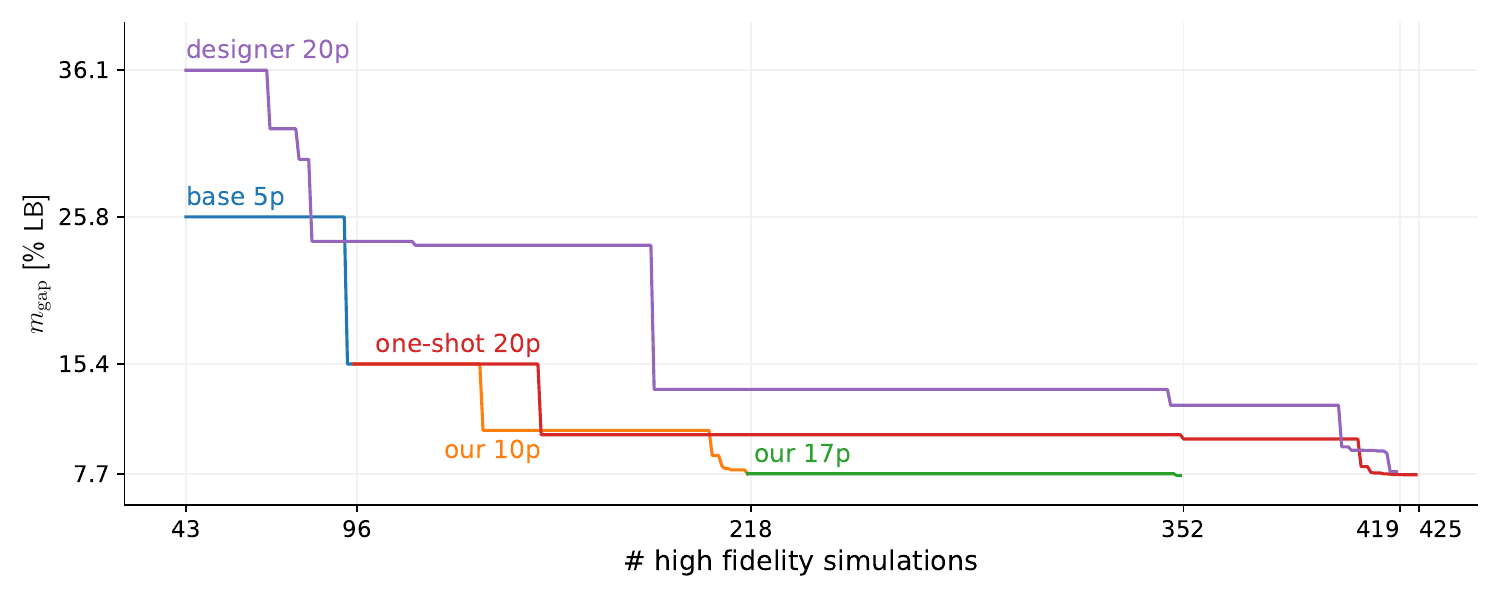}
	\caption{Evolution of the optimal high-fidelity configuration for different parameterizations of
the midship section.}
	\label{fig:msection_hifi_mass_history}
\end{figure}

The evolution of optimal high-fidelity configurations for all the approaches is shown in  Figure~\ref{fig:msection_hifi_mass_history}.
The large sequences of non-decreasing ships are due to the multi-objective optimization, in which the penalized mass is only incidentally optimized and the focus is on obtaining better exploration of the PF.
Table~\ref{tab:msection_results} collects the results of the optimum for each problem formulation.
Overall, all the largest parameterizations achieve similar results in terms of gap from the lower bound, but with a number of differences.
Iterative reparameterization shows a much faster decrease, with the $10$ parameters model reaching a comparable objective value to the others, but with a much lower number of high-fidelity simulations.
The second reparameterization, with $17$ parameters, decreases the optimum only by a negligible amount.
However, we observe that the number of buckled elements is about $11\%$ lower than the threshold value, suggesting that further rework by a human designer could benefit from this buffer from the constraint's critical value.
The one-shot parameterization requires a larger number of high-fidelity simulations in the multi-objective phase, but the buffer in terms of buckled elements is even more pronounced than for the $17$ parameter model, being about $32\%$ of the threshold value.
The designer-provided reference requires a number of high-fidelity evaluations similar to the one-shot case, but the final configuration reaches the threshold on the number of buckled elements.

\begin{table}[hbt!]
	\caption{Optimization results for different parameterizations of the midship section. The thresholds for the number of yielded and buckled elements are 200 and 4000, respectively. The threshold for the deflection is \SI{140}{mm} and the maximum VCG is \SI{15}{m}.}
    \label{tab:msection_results}
	\centering
	\begin{tabular}{ l r r r c c }
        \hline
        \hline
        Parameterization & $m_{\text{gap}}$ [\% LB] & \#yielded & \#buckled & Deflection [\si{mm}] & VCG [\si{m}] \\
        \hline
        \hline
        \rowcolor{Gray}
        Initial configuration       & $>1000$ & 395 & 6211 & 150 & 14.79 \\
        Base 5 p.                   & 15.44   & 196 & 3994 & 138 & 14.93 \\
        \rowcolor{Gray}
        Our method 10 p.            & 7.70    & 33  & 3984 & 140 & 14.99 \\
        Our method 17 p.            & 7.57    & 11  & 3563 & 140 & 14.99 \\
        \rowcolor{Gray}
        Our one-shot 20 p.          & 7.64    & 5   & 2729 & 140 & 14.99 \\
        Designer 20 p.              & 7.86    & 22  & 4004 & 140 & 14.99 \\
        \hline
        \hline
	\end{tabular}
\end{table}

\subsection{Full ship}
\label{subsec:results_fullship}
\begin{figure}[hbt!]
	\centering
	\includegraphics[trim=0 0 0 0, clip, width=.48\textwidth]{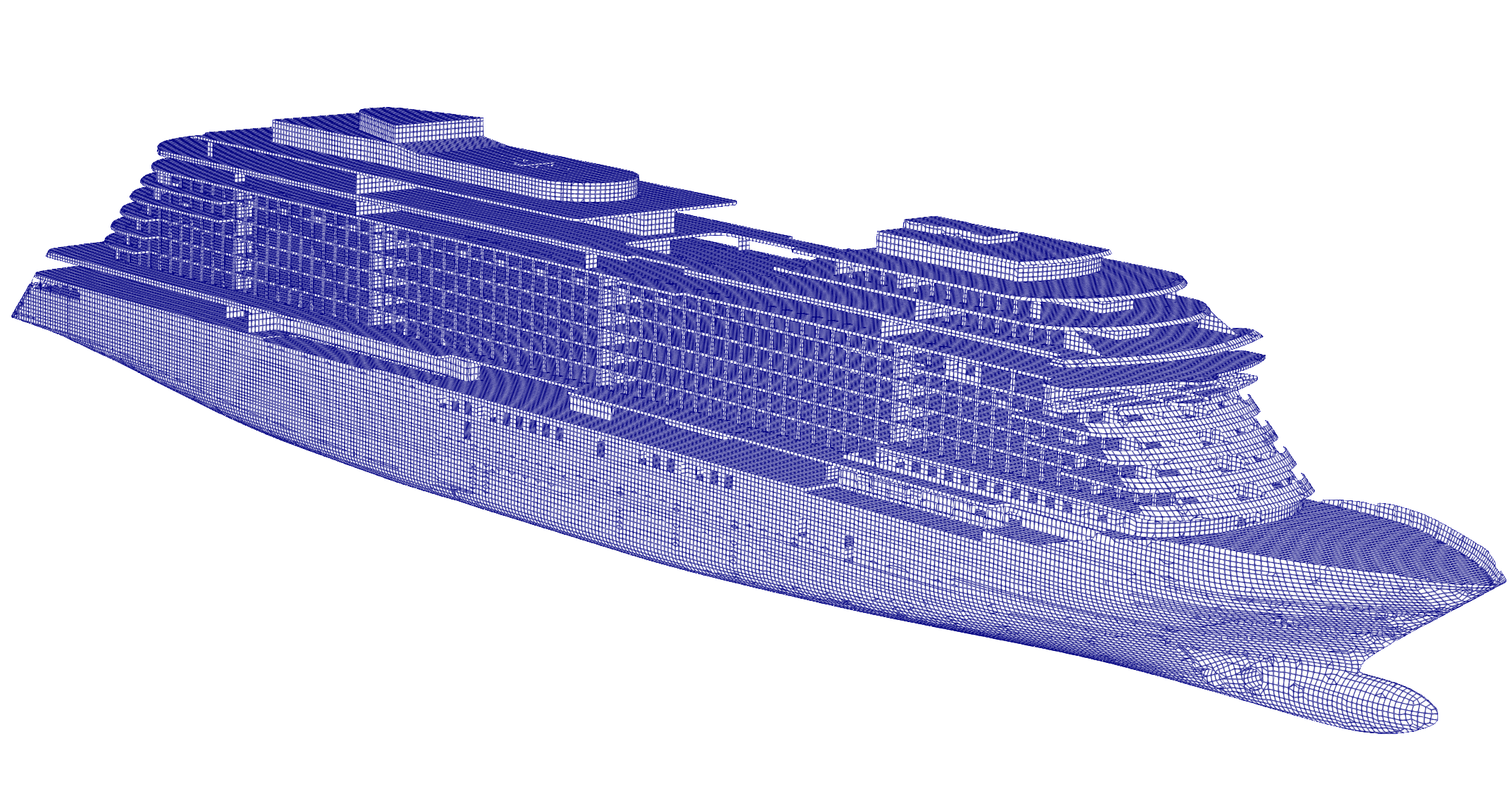}\hfill
	\includegraphics[trim=0 0 0 0, clip, width=.48\textwidth]{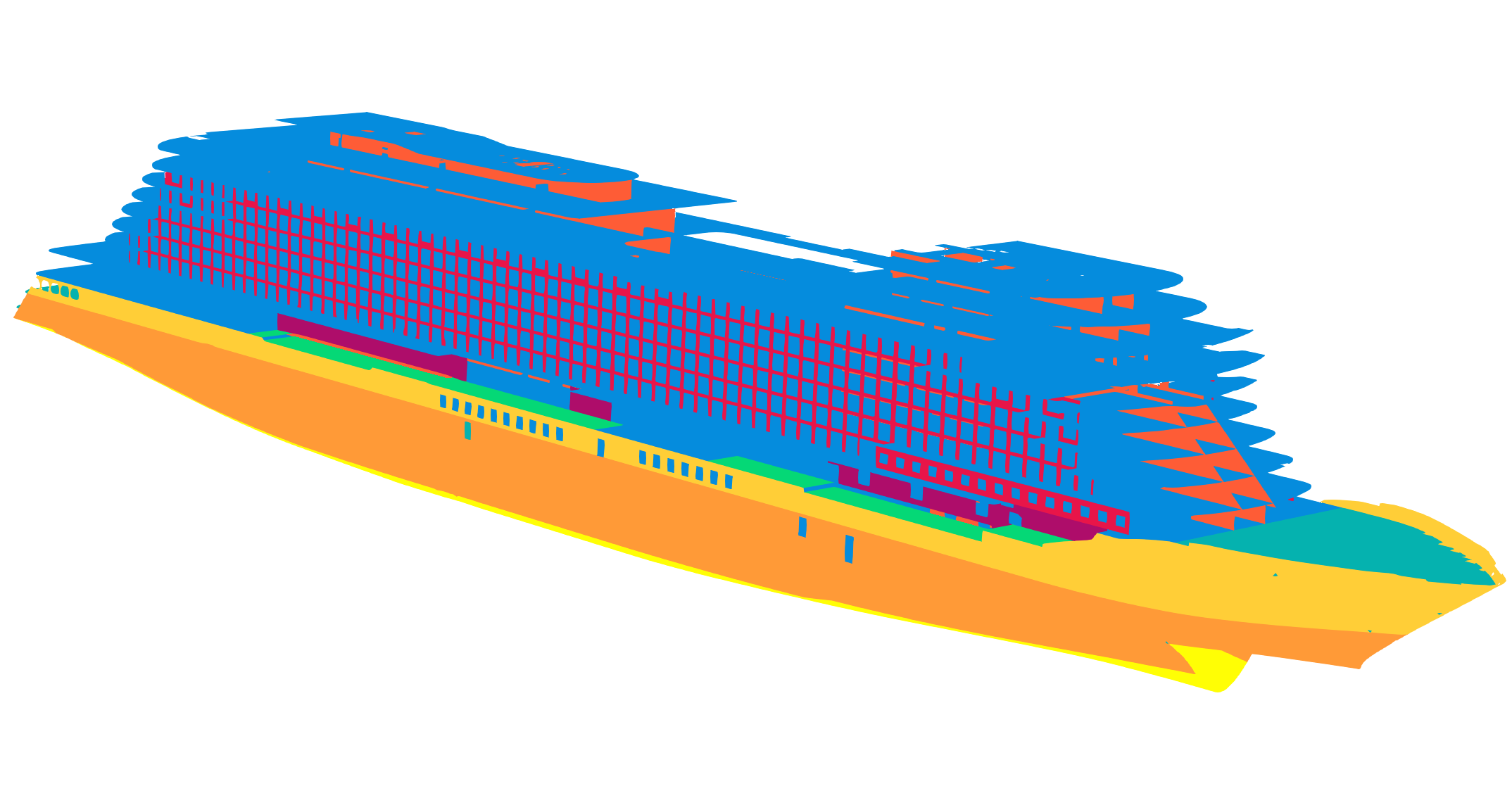}
    \caption{Full view of the ship on the left, initial parameterization on the right.}
	\label{fig:fullship_general}
\end{figure}

The full scale test model in Figure~\ref{fig:fullship_general} comes from a ship developed by Fincantieri SpA. 
This model consists of $485736$ elements, of which, $270656$ are linked to the parameters in Table~\ref{tab:fullship_params} and organized in $5005$ patches. 
The MSC NASTRAN solution time for this model is $7$ minutes, and the surrogate prediction time is $0.5$ seconds, for a speedup factor of $840$. 
The default parameter configuration corresponds to the minimum regulatory thickness of the parameterized sections, while all the upper bounds have been left to the maximum commercially available value. 
The parameters are equally divided into two groups, one containing decks and bulkheads having a high number of feasible values, and the rest of the structural elements with small domain. The resulting domain contains 467 million unique combinations. %
Due to the increased complexity of this model, the maximum number of parameters is $40$.
To comply with the proprietary constraints set by the industrial partner, we present the QoIs in a relative format and the critical values are not shown.

\begin{table}[htb]
	\caption{Parameters description of the full ship test case. Thickness values are in $\si{mm}$.}
    \label{tab:fullship_params}
	\centering
	\begin{tabular}{ c l r r r r }
		\hline
		\hline
		Parameter & Region & Patches & Min. & Max. & \# values \\
		\hline
		\hline
		\rowcolor{Gray}
		$\xpar_{1}$ & Decks in public areas    & 3038 & 5.0  & 25.0 & 16 \\
		$\xpar_{2}$ & Decks in machinery areas & 548  & 6.5  & 25.0 & 13 \\
		\rowcolor{Gray}
		$\xpar_{3}$ & Lifeboat exposed deck    & 90   & 13.0 & 25.0 & 5 \\
		$\xpar_{4}$ & Inner bottom             & 252  & 12.0 & 25.0 & 6 \\
		\rowcolor{Gray}
		$\xpar_{5}$ & Bottom                   & 312  & 15.0 & 25.0 & 3 \\
		$\xpar_{6}$ & Shell plating, above waterline  
                                               & 60   & 14.0 & 25.0 & 4 \\
		\rowcolor{Gray}
		$\xpar_{7}$ & Shell plating, below waterline 
                                               & 118  & 15.0 & 25.0 & 3 \\
		$\xpar_{8}$ & Internal bulkheads       & 367  & 5.0  & 25.0 & 16 \\
		\rowcolor{Gray}
		$\xpar_{9}$ & External bulkheads, in super-structures 
                                               & 190  & 8.0  & 25.0 & 10 \\
		$\xpar_{10}$ & External bulkheads      & 31   & 6.5  & 25.0 & 13 \\
		\hline
		\hline
	\end{tabular}
\end{table}

The initial sampling of the parameter space is limited to $20$ configurations other than the default one. 
The truncation rank chosen for the POD is 16, and the GPRs use the same structure as for the midship section case.
Details on the truncation ranks and the corresponding errors can be found in Appendix~\ref{subsec:appendix_fullship}.
As before, multi-objective optimization is carried out in order to identify the appropriate thresholds for the number of yielded and buckled elements. 
The maximum value of VCG is determined by the model geometry, while no constraint on the vertical deflection is specified.
We repeat the same tests from Section~\ref{subsec:results_msection}: a sequence of iterative reparameterizations from $10$ to $40$ parameters, a one-shot reparameterization from $10$ to $40$ parameters, and the optimization of a $40$ parameters model provided by the designers.

\begin{figure}[hbt!]
	\centering
	\includegraphics[trim=0 0 0 0, clip, width=\textwidth]{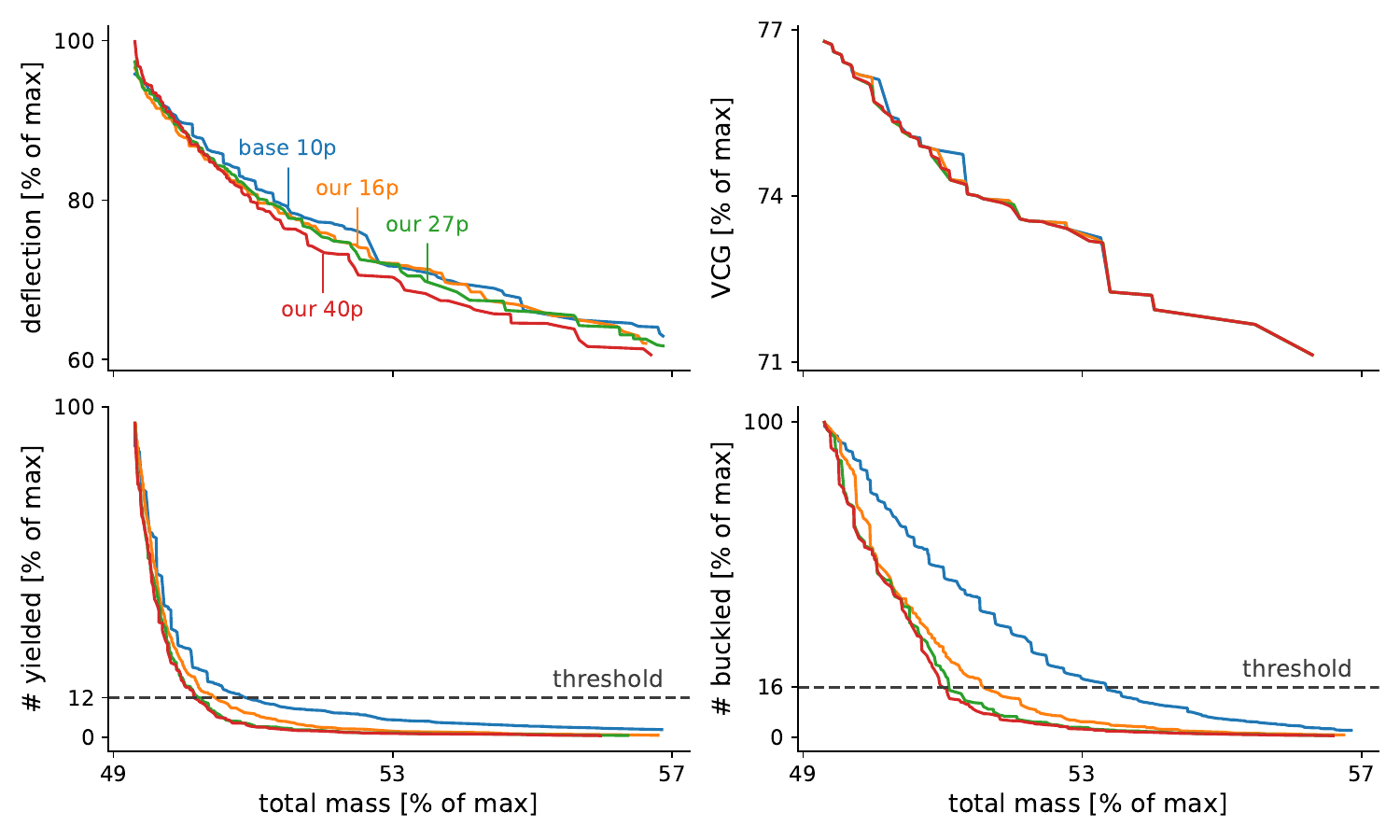}
	\caption{Evolution of the surrogate Pareto frontier for the full ship, for different parameterizations.
    The threshold for the VCG constraint is not shown, since the PF of the projection lies well below it.}
	\label{fig:fullship_pareto_frontiers}
\end{figure}

In the iterative reparameterization test, the refinement step is performed $3$ times giving models with $16$, $27$ and $40$ parameters each.
Figure~\ref{fig:fullship_pareto_frontiers} shows the evolution of the surrogate PF.
Each refined model improves the previous PF, with the largest effect being achieved on the number of buckled elements.
The first reparameterization shows a more pronounced reduction of the objective, while the third appears to have a limited effect.
The threshold values for the number of yielded and buckled elements result much more stringent than those used in Section~\ref{subsec:results_msection}.

\begin{figure}[hbt!]
	\centering
	\includegraphics[trim=0 0 0 0, clip, width=\textwidth]{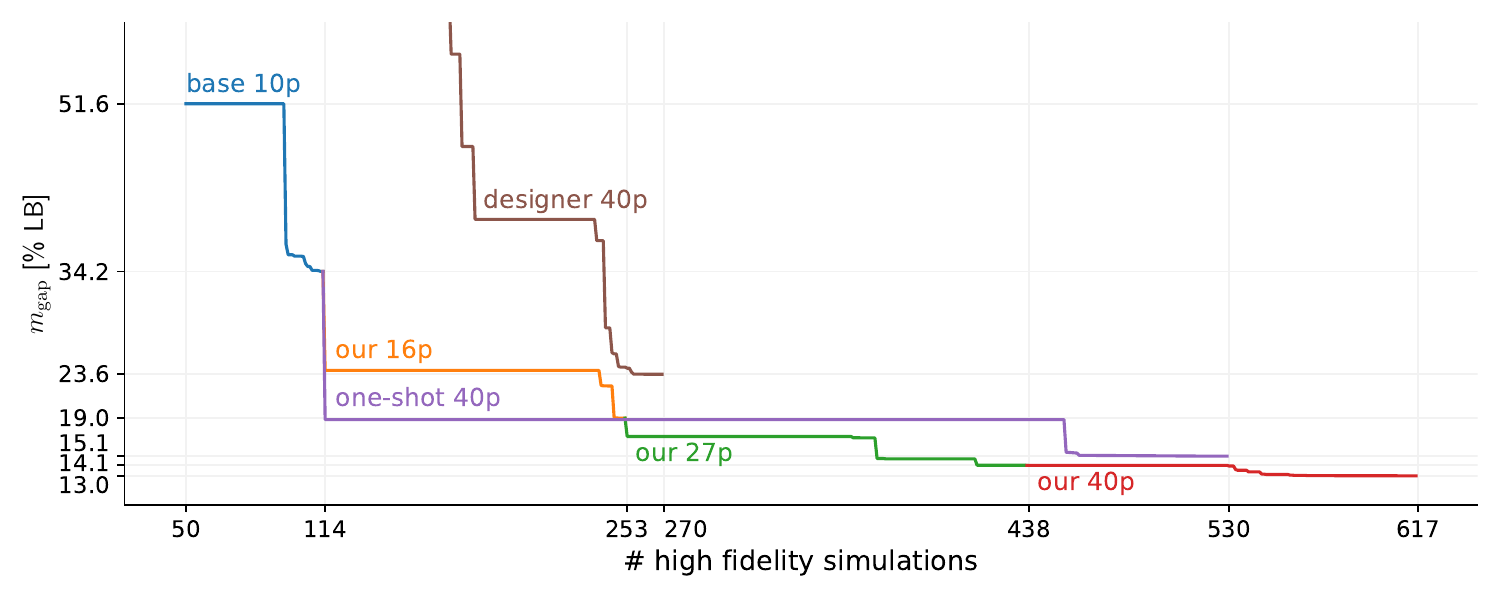}
	\caption{Evolution of the optimal high-fidelity configuration for different parameterizations of the full ship.}
	\label{fig:fullship_hifi_mass_history}
\end{figure}

The evolution of optimal high-fidelity configurations for all the approaches is shown in Figure~\ref{fig:fullship_hifi_mass_history}, and the details of the optimum configuration for each problem is shown in Table~\ref{tab:fullship_results}.
A striking difference from Figure~\ref{fig:msection_hifi_mass_history} is the large decrease in objective function due to the reparameterizations to $16$ and $40$ parameters, showing that the ILP was able to generate a highly performant and feasible configuration.
As in the case of the PFs, refinements beyond the first incur in lower gains, and the final refinement of the incremental approach only decreases the optimum value by around $1\%$ of the LB defined in Equation~\eqref{eq:mass_gap}.
Multi and single-objective optimization of the one-shot model are not able to achieve as large improvements as the initial ILP-optimal configuration, with the final result being worse than the $27$ parameters model by about $1\%$ of the LB.
The final configuration of the designer model is better than the $10$ parameter model and reaches its optimum with a lower number of high-fidelity simulations compared with the other approaches, but it is only marginally better than the ILP optimum that generated the $16$ parameters model.
All the parameterizations generated through our automated procedure were able to outperform the one from the designers.
These results suggest that the high complexity of the full ship, with multiple large parameterized groups of patches, provided a harder challenge in designing an effective parameterization, even for an expert.

Regarding the constrained QoIs, the number of buckled elements and VCG are close to the threshold for all instances.
The number of yielded elements does so only for the models that were iteratively refined to $27$ and $40$ parameters, while in the other instances it remains around half of the threshold value.
This behavior is different from what was observed in Table~\ref{tab:msection_results}, where both QoIs related to element failure were consistently lower than the threshold for the models above 10 parameters.
Since the incremental approach is able to produce models with lower objective function and higher activation of the constraints, we conclude that it is the most suited for the initial design phase.

\begin{table}[hbt!]
	\caption{Optimization results for different parameterizations of the full ship.}
    \label{tab:fullship_results}
	\centering
	\begin{tabular}{ l r r r r }
        \hline
        \hline
        Parameterization & $m_{\text{gap}}$ [\% LB] & \# yielded [\% crit] & \# buckled [\% crit] & VCG [\% crit] \\
        \hline
        \hline
        \rowcolor{Gray}
        Initial configuration       & $>1000$   & 789.75 & 636.38 & 96.26 \\
        Base 10 p.                  & 34.21     & 46.25  & 100.38 & 99.35 \\
        \rowcolor{Gray}
        Our method 16 p.            & 19.01     & 56.00  & 99.88  & 99.76 \\
        Our method 27 p.            & 14.12     & 98.25  & 99.29  & 99.71 \\
        \rowcolor{Gray}
        Our method 40 p.            & 13.03     & 98.25  & 100.01 & 99.88 \\
        Our one-shot 40 p.          & 15.07     & 51.50  & 99.87  & 99.85 \\
        \rowcolor{Gray}
        Designer 40 p.              & 23.55     & 55.00  & 100.15 & 99.99 \\
        \hline
        \hline
	\end{tabular}
\end{table}

\section{Conclusions and perspectives}
\label{sec:conclusions}
In this paper, we presented a new hierarchical parameterization refinement method that leverages ROMs and a series of optimization subproblems. The parameterization does not rely on specific optimization methods and can be adapted to other use cases. We showcased the benefit of such a technique by extending an existing pipeline for the automated structural optimization of cruise ships with the addition of multi-objective optimization capabilities, a specialized BO procedure, a local search heuristic, and the aforementioned parameterization refinement procedure.

The multi-objective optimization is based on the established NSGA genetic algorithm, coupled with an infill criterion leveraging the GPR component of the surrogates.
The Bayesian optimization procedure has been enhanced with the addition of linear inequality constraints on the incumbent and the VCG, and an ILP-based rounding procedure to overcome the rounding issues due to the discrete nature of the parametric space.
The principal dimensions search provides a local, black-box refinement procedure to further optimize the surrogates when BO is not able to find new candidates within the computation budget.
Finally, the reparameterization module refines the optimization problem formulation by adapting to the emergent behavior of the QoIs.
New parameters are constructed by finding a clustering of patches in each parameterized section, to separate the responses to thickness changes in terms of yielded and buckled elements.
The clusters are determined through the solution of ILPs constructed on the surrogate models' predictions.

The hierarchical nature of the refinement enables efficient use of past expensive high-fidelity results.
The enhanced pipeline has been tested on the mass optimization of two test cases, a midship section and a large-scale full ship, for which an initial parameterization was refined both iteratively and at once, and the final results were compared with a baseline given by expert designers.
In the simpler case, the pipeline results are comparable to the baseline even though a lower number of parameters was used, suggesting the effectiveness of the approach.
In the more complex case, the pipeline was able to achieve a much lower total ship mass than the baseline, with a large gain due to the configurations found by the clustering procedure.
The proposed pipeline has proven effective in achieving a ship configuration with a low total mass, while respecting the safety and manufacturing constraints imposed by the designers.
The use of this pipeline can help streamline the initial design phase of complex or novel ships, providing good configurations with minimal supervision.

Further directions of development could be the integration of other types of decision variables, such as steel quality, or the addition of different buckling models to support more classification rules.
Another extension could be the introduction of fatigue modeling, which would however require adaptation of the surrogate construction and evaluation.
Although formulated for the optimization of cruise ships, the pipeline could be applied to the structural design phase of other structures where the limiting factor is the validation of static yielding and buckling phenomena.
QoIs such as the vertical deflection and VCG could be substituted by other quantities relevant to the new context.
The reparameterization module, however, relies only on the type of assignment of the decision variables to the geometrical elements and can be applied without changes to new problems of the same nature.

\section*{Acknowledgements}
This work was partially supported by the PNRR (PhD with industries) project "Development of advanced reduced methods of parametric optimization for the structural analysis of cruise vessels", co-sponsored by Fincantieri SpA, by the PNRR NGE project iNEST "Interconnected Nord-Est Innovation Ecosystem", and by the European Research Council Executive Agency by the Proof of Concept project ARGOS, "Advanced Reduced Groupware Online Simulation", P.I. Professor Gianluigi Rozza.

\bibliographystyle{abbrvurl}
\bibliography{main_arxiv}

\appendix

\section{Evaluation of the errors in the surrogate models}
Surrogate models are used to quickly select promising configurations for the minimization of the QoIs.
Due to the expensive validation of the results using the high-fidelity solver, the database of configurations, from which the ROMs are constructed, mostly comprises the results from multi- and single-objective optimizations, as selected by their infill criteria.
The only exceptions come from the initial random sampling, and from the reparameterization steps, where the larger domain is again randomly sampled before progressing to the optimization.

We evaluate the choice of the truncation rank used in the POD, and its update after the reparameterization.
Automatization of these steps will make the pipeline more user-friendly and provide smoother integration with existing workflows.
We report on the behavior of the POD singular values and the prediction errors for different combinations of truncation rank and high-fidelity database sizes.
In Section~\ref{subsec:appendix_midship} we analyze the models based on the midship section, while Section~\ref{subsec:appendix_fullship} presents the full ship results.

\subsection{Midship section}
\label{subsec:appendix_midship}
\begin{figure}[hbt!]
	\centering
	\includegraphics[trim=0 0 0 0, clip, width=1\textwidth]{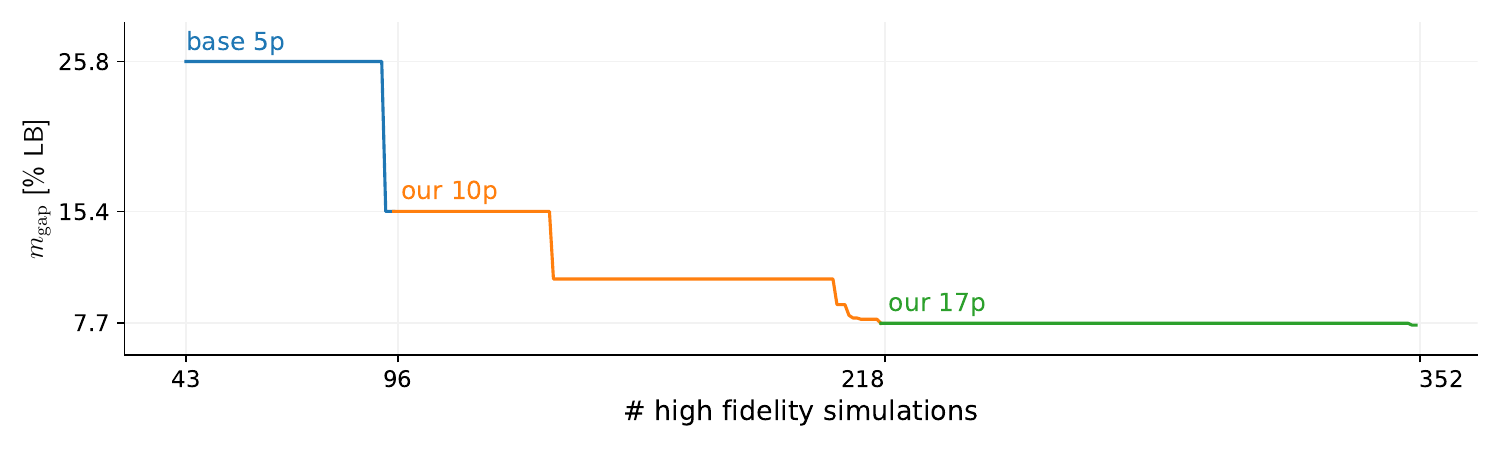}
	\caption{Evolution of the high fidelity optimal configuration for the midship section, under the iterative reparameterization approach. For a comparison with different settings see Figure~\ref{fig:msection_hifi_mass_history}.
    }
	\label{fig:msection_hifi_mass_history_redux}
\end{figure}

Figure~\ref{fig:msection_hifi_mass_history_redux} shows the evolution of the high-fidelity optimum for the midship section, for the iterative reparameterization approach described in Section~\ref{subsec:results_msection}. 
The plot begins at the 43rd sample, as the earlier ones are significantly higher than the final optimum.
Multi-objective optimization provides the largest component of the high-fidelity database, with configurations spread throughout the domain, while single-objective optimization limits configurations close to the latest optimum.

\begin{figure}[hbt!]
	\centering
	\includegraphics[trim=0 0 0 0, clip, width=.93\textwidth]{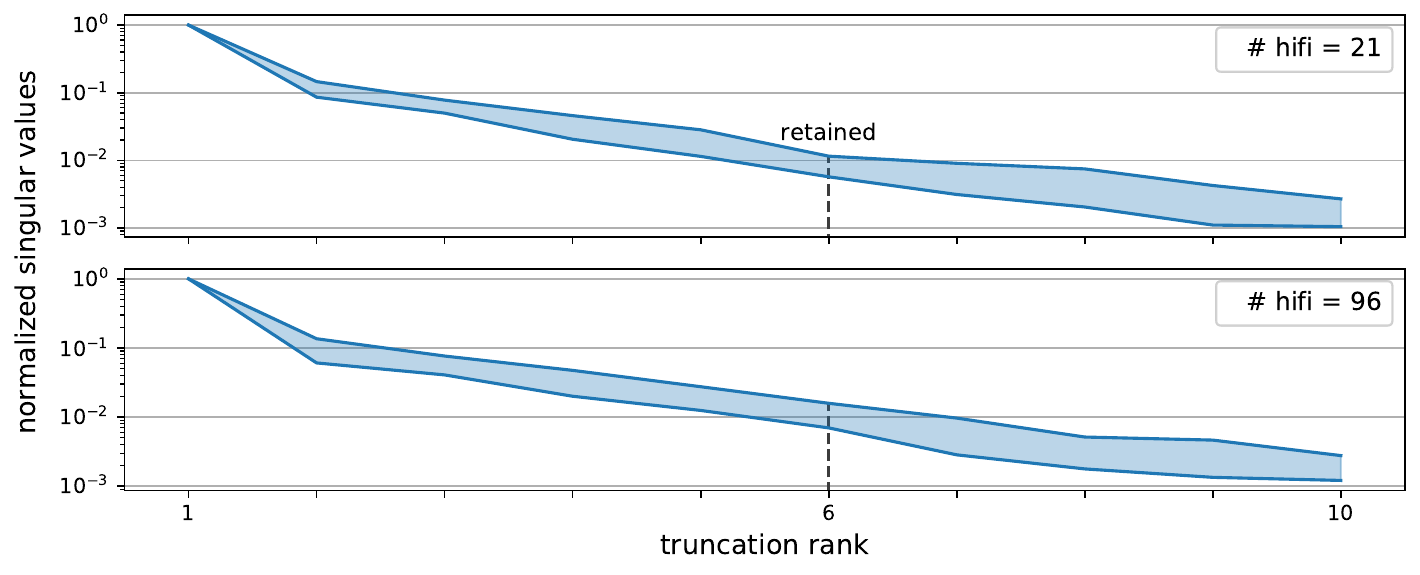}
	\caption{Evolution of the minimum and maximum for the normalized singular values decay of the stress tensor components, for the midship section with 5 parameters.
    Values are reported for the high-fidelity database before and after optimization.
    }
	\label{fig:msection_5p_svalues_evo}
\end{figure}

The choice of truncation rank for the POD happens after the initial random sampling.
Figure~\ref{fig:msection_5p_svalues_evo} shows the evolution of the singular values decay for the stress tensor components, for different database sizes.
The databases considered are the initial sampling (the initial configuration and 20 random samples) and the final one, after single-objective optimization.
For the initial sampling, the normalized singular values from rank 7 and above fall below the $10^{-2}$ threshold and are thus discarded.
As new high-fidelity samples are acquired, the largest $7$-th singular value never exceeds the $10^{-2}$ threshold, confirming that the truncation choice was robust.

\begin{figure}[hbt!]
	\centering
	\includegraphics[trim=0 0 0 0, clip, width=.97\textwidth]{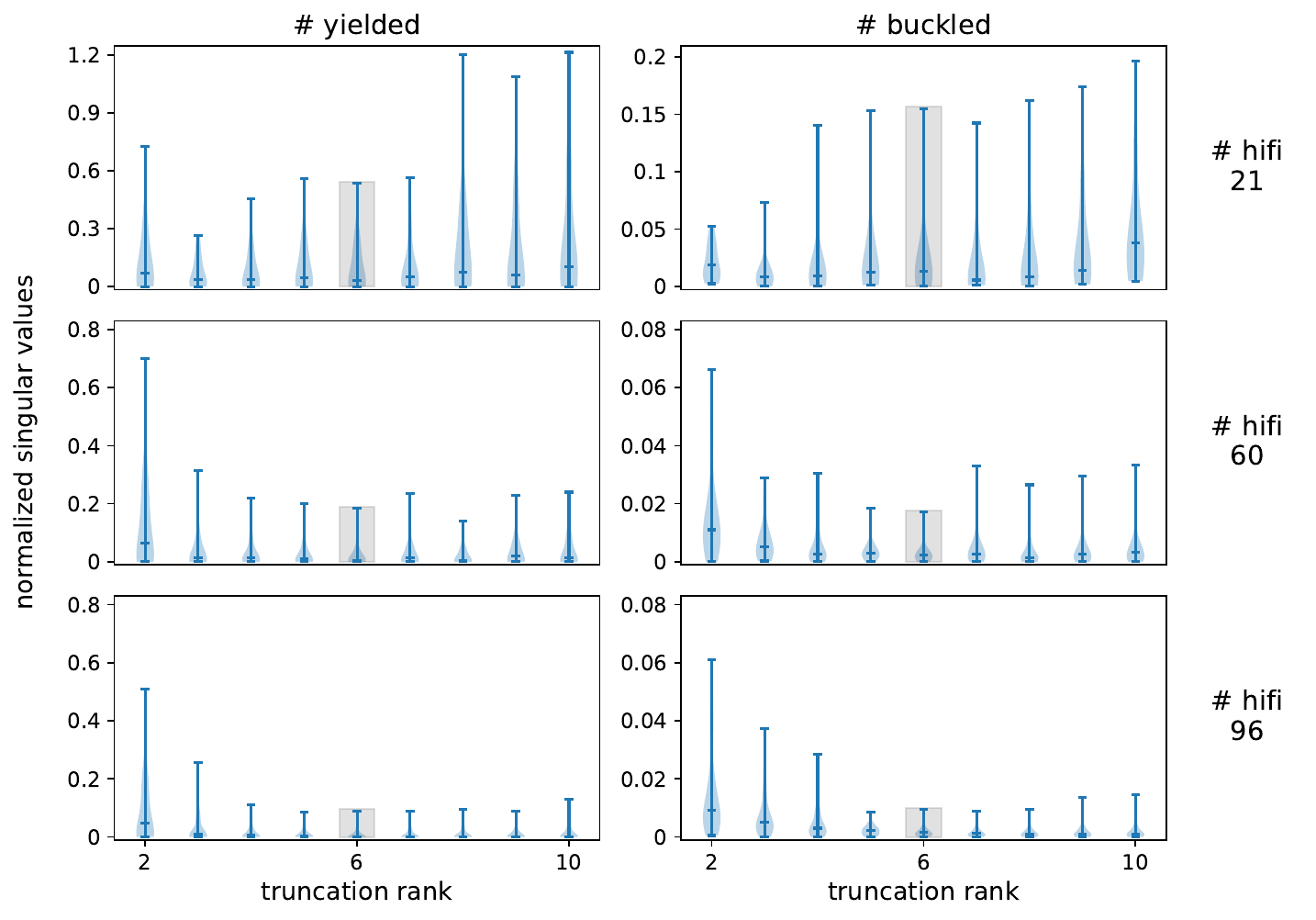}
	\caption{Evolution of the 5-fold cross-validation error distribution for the midship section with 5 parameters, for different sizes of the high-fidelity database. 
    For the number of yielded and buckled elements, the difference between the predicted and true values is divided by the critical threshold. 
    The plots highlight the median error and the gray box identifies the chosen truncation rank.
    }
	\label{fig:msection_5p_errors_evo}
\end{figure}

Figure~\ref{fig:msection_5p_errors_evo} shows the behavior of the 5-fold cross-validation error on the number of yielded and buckled elements, as a fraction of the critical threshold (200 for the number of yielded elements, 4000 for the buckled elements).
We also show the results for a database approximately in the middle of the procedure, during the multi-objective optimization phase.
The initial database does not show a clear decrease of prediction errors as the truncation rank grows, while during the middle of the pipeline, the surrogates achieve good predictive performances.
With the largest database, the surrogates show little difference between truncation ranks 6 and 7.
The choice of truncation rank 6 appears to be appropriate for the initial parameterization of the midship section.

\begin{figure}[hbt!]
	\centering
	\includegraphics[trim=0 0 0 0, clip, width=.93\textwidth]{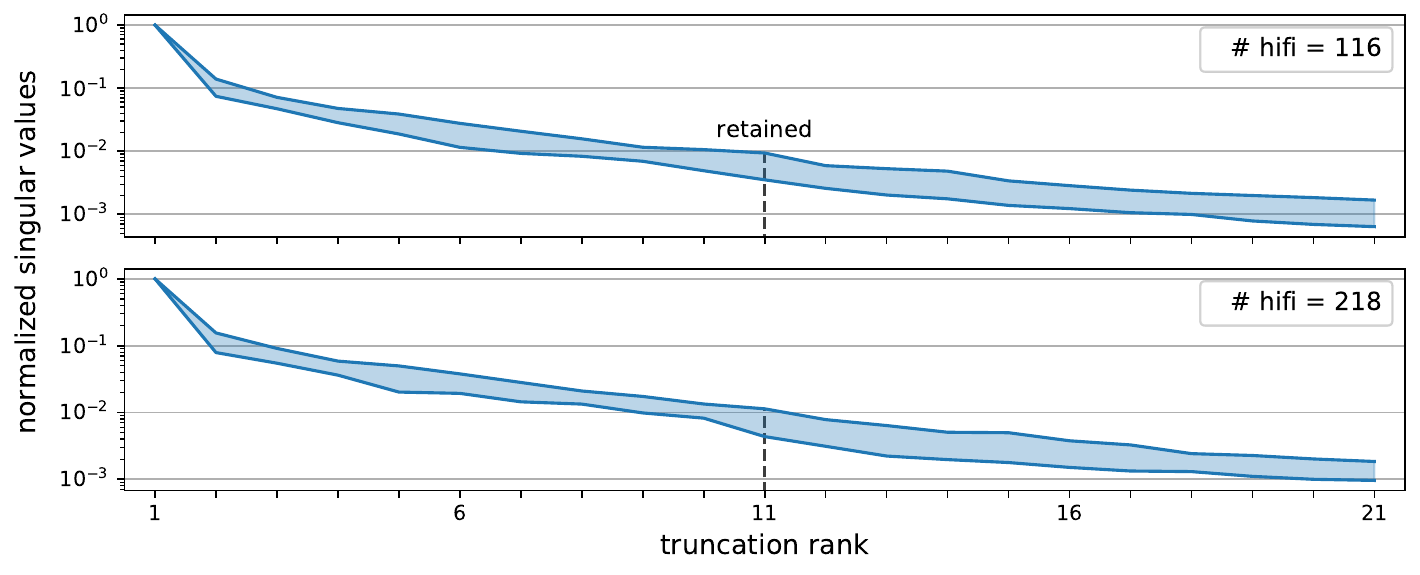}
	\caption{Evolution of the minimum and maximum for the normalized singular values decay of the stress tensor components, for the midship section with 10 parameters.
    Values are reported for the high-fidelity database before and after optimization.
    }
	\label{fig:msection_5p_to_10p_svalues_evo}
\end{figure}

\begin{figure}[hbt!]
	\centering
	\includegraphics[trim=0 0 0 0, clip, width=.97\textwidth]{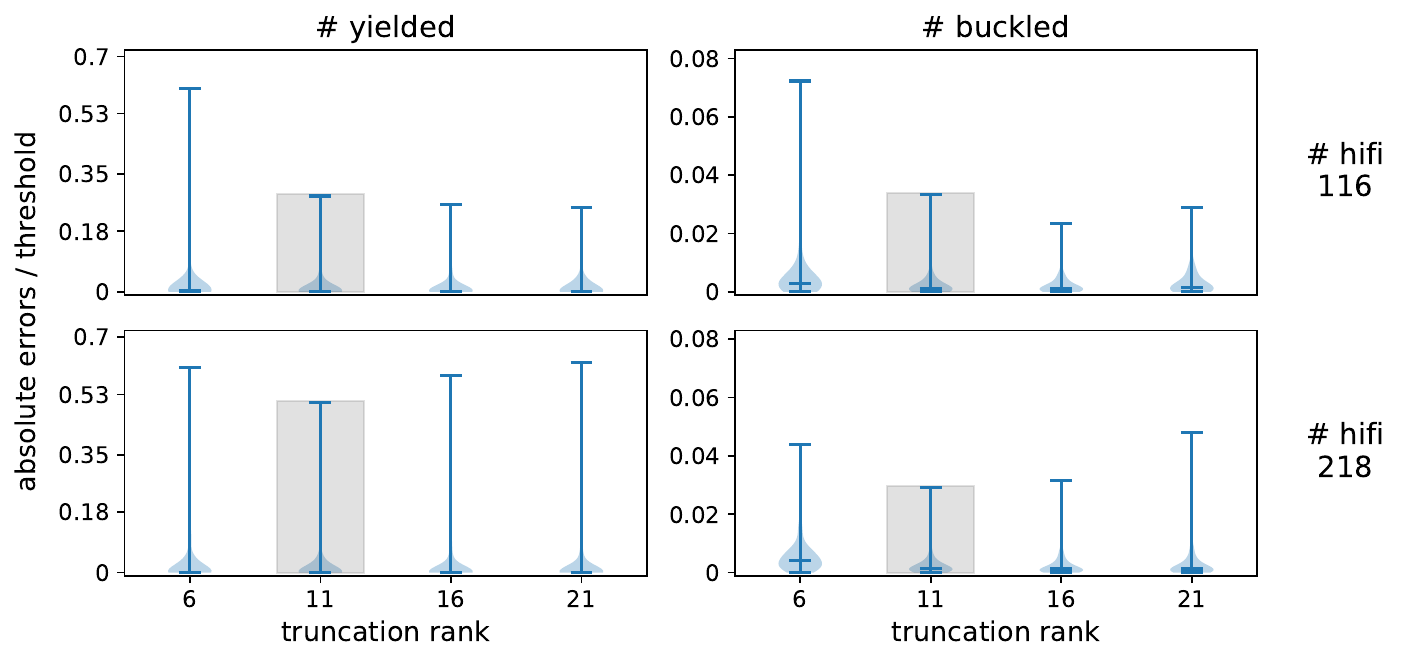}
	\caption{Evolution of the 5-fold cross-validation error distribution for the midship section with 10 parameters.
    Values are reported for the high-fidelity database before and after optimization. The gray box identifies the chosen truncation rank.
    }
	\label{fig:msection_5p_to_10p_errors_evo}
\end{figure}

Figures~\ref{fig:msection_5p_to_10p_svalues_evo} and \ref{fig:msection_5p_to_10p_errors_evo} show the same information discussed before, but for the model after its reparameterization to 10 parameters.
The first database considered contains the previous 96 samples, and 20 random samples from the enlarged domain.
The truncation rank was automatically increased from 6 to 11, and the first POD of the stress tensors shows how the largest normalized singular values reach the $10^{-2}$ threshold between ranks 9 and 11.
Further enrichments of the high-fidelity database determine a slower decay, but the truncation rank at 11 still discards the modes with a normalized singular value below $10^{-2}$.
As for the prediction error of yielded and buckled elements, the first surrogates constructed on the enlarged domain show an increase in the maximum error values compared to the latest surrogates using 5 parameters.
The truncation rank 11 is compared to 6, 16, and 21, showing that retaining a larger number of modes does not provide substantial gains in accuracy, and suggesting instead an increase in prediction errors due to outliers.
Performances are comparable to the 5-parameters models in term of median values.

\begin{figure}[hbt!]
	\centering
	\includegraphics[trim=0 0 0 0, clip, width=.93\textwidth]{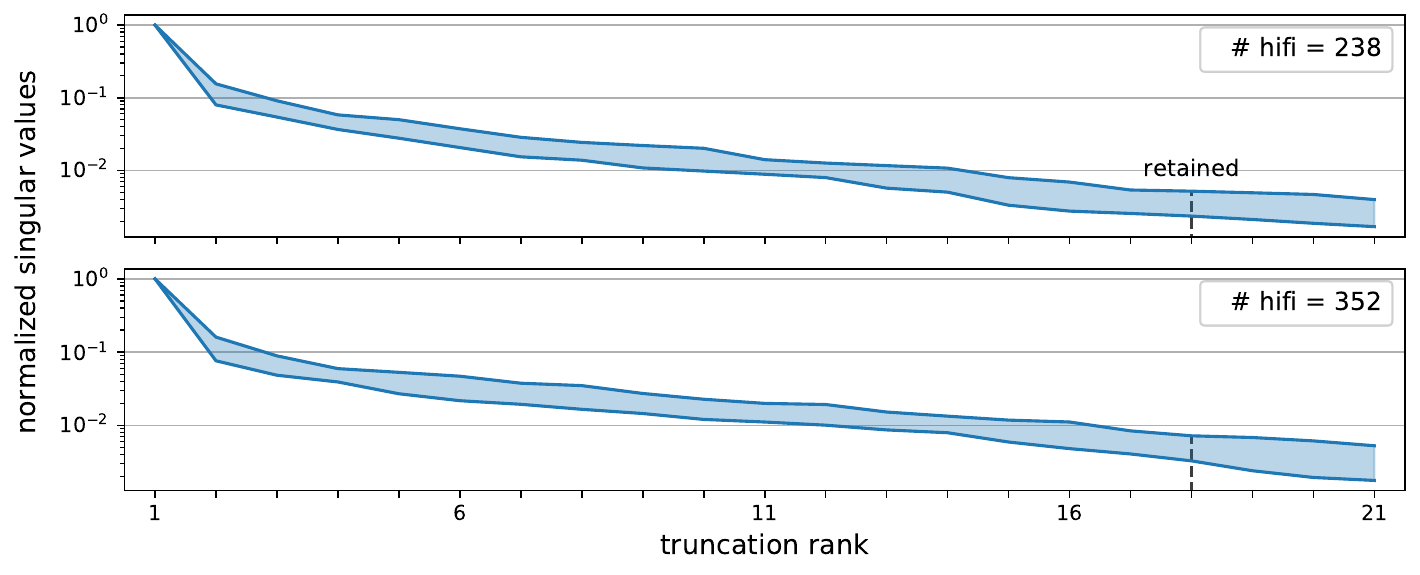}
	\caption{Evolution of the minimum and maximum for the normalized singular values decay of the stress tensor components, for the midship section with 17 parameters.
    Values are reported for the high-fidelity database before and after optimization.
    }
	\label{fig:msection_10p_to_17p_svalues_evo}
\end{figure}

\begin{figure}[hbt!]
	\centering
	\includegraphics[trim=0 0 0 0, clip, width=.97\textwidth]{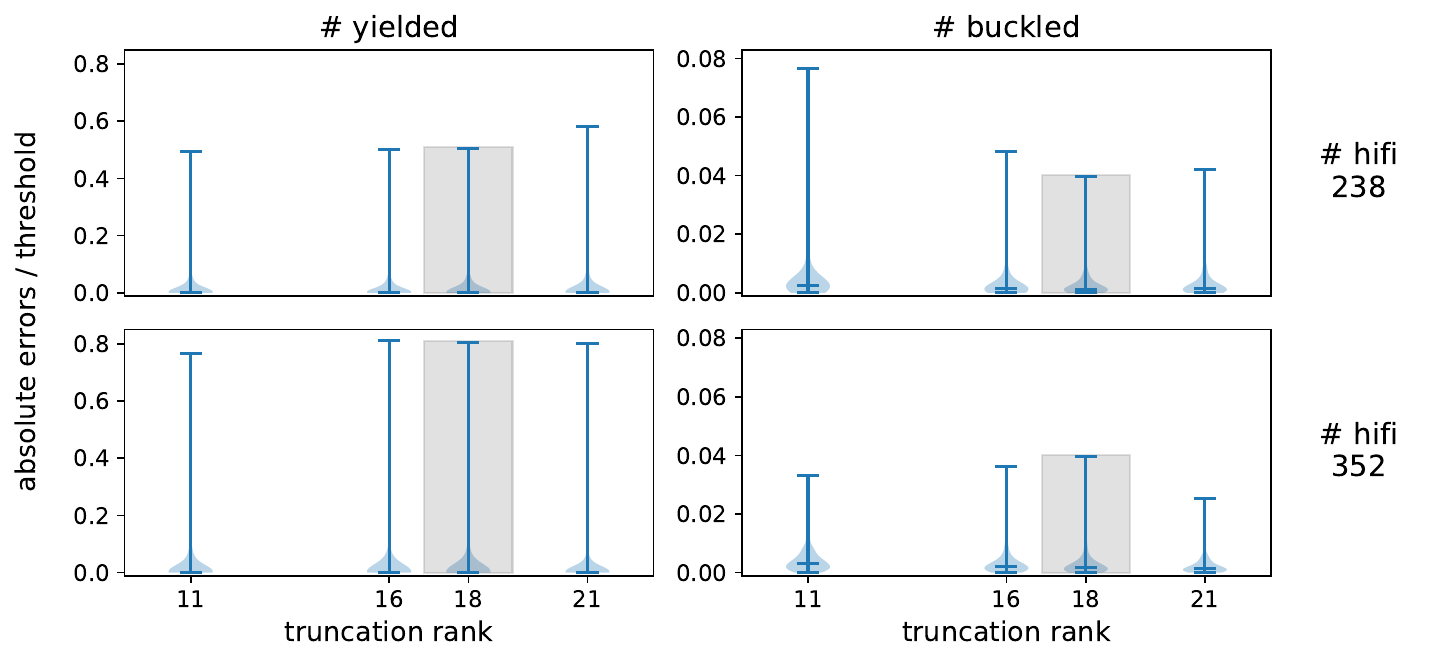}
	\caption{Evolution of the 5-fold cross-validation error distribution for the midship section with 17 parameters.
    Values are reported for the high-fidelity database before and after optimization. The gray box identifies the chosen truncation rank.
    }
	\label{fig:msection_10p_to_17p_errors_evo}
\end{figure}

With the second reparameterization of the model, the new domain has 17 dimensions and the truncation rank is increased to 18.
Again, the initial database contains the 218 configurations from the previous optimization, and 20 samples from the new domain.
This time, the normalized singular values decay in Figure~\ref{fig:msection_10p_to_17p_svalues_evo} crosses the $10^{-2}$ threshold between ranks 14 and 15, and further exploration of the domain only moves it between ranks 16 and 17.
In this case, the automatic update of the truncation rank was conservative.
The behavior of the error on the number of yielded elements does not show a marked dependence on the truncation ranks tested, as shown in Figure~\ref{fig:msection_10p_to_17p_errors_evo}.
As for the error on the number of buckled elements, the largest high-fidelity database shows only a marginally lower median error from rank 16 and above, compared to rank 11.

The results from the evaluation validate the choice of $10^{-2}$ as the threshold for the truncation rank in the POD, at the start of the optimization pipeline.
This value is unusually high compared to the common practice, but as shown by the results for the prediction errors, choosing a higher truncation rank seems to decrease, rather than increase, the accuracy of the ROMs for evaluating the QoIs.
A possible explanation is due to the combination of an increase in the complexity of performing regression on a larger latent space, and the non-linearity of the post-processing step for computing the QoIs, especially the buckling state of the elements.
Indeed, the singular values decay does not provide information on the variability of the coefficients across the domain nor in the areas close to the best-performing configurations, so reconstruction errors in latent space could lead to larger effects on the QoIs.

\subsection{Full ship}
\label{subsec:appendix_fullship}
\begin{figure}[ht!]
	\centering
	\includegraphics[trim=0 0 0 0, clip, width=1\textwidth]{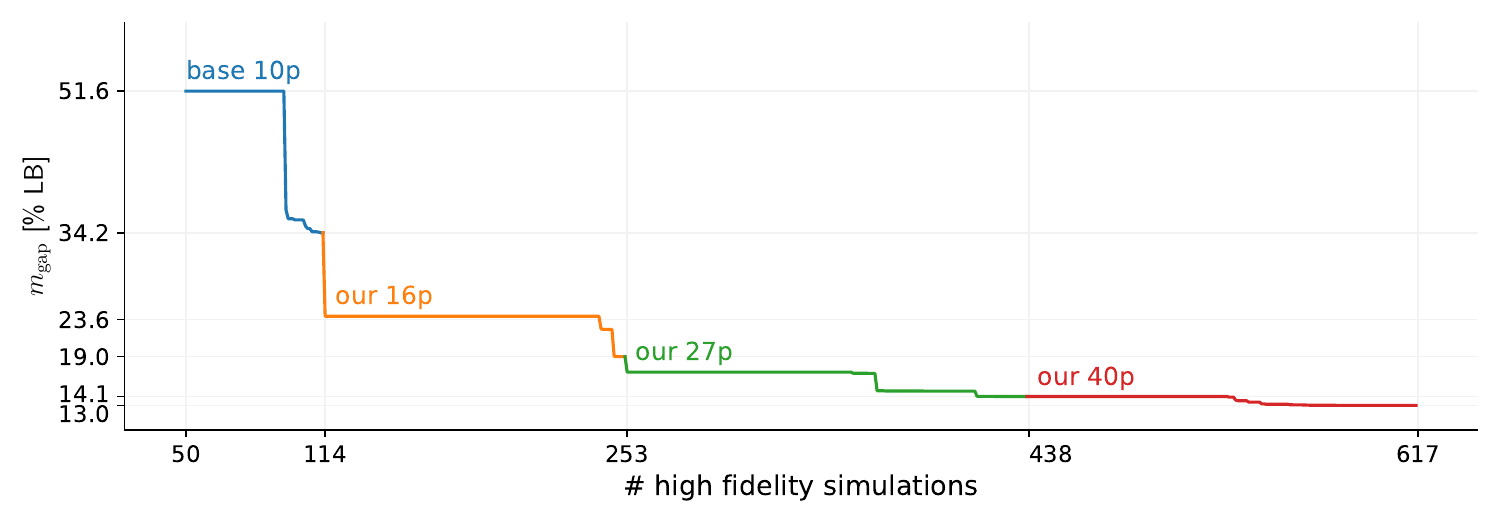}
	\caption{Evolution of the high fidelity optimal configuration for the full ship, under the iterative reparameterization approach. For a comparison with different settings see Figure~\ref{fig:fullship_hifi_mass_history}.
    }
	\label{fig:fullship_hifi_mass_history_redux}
\end{figure}

The evolution of the high-fidelity optimum, during the iterative reparameterization of the full ship model, is depicted in Figure~\ref{fig:fullship_hifi_mass_history_redux}.
The plot begins at the 50th sample, as the earlier ones are significantly higher than the final optimum.

\begin{figure}[hbt!]
	\centering
	\includegraphics[trim=0 0 0 0, clip, width=.93\textwidth]{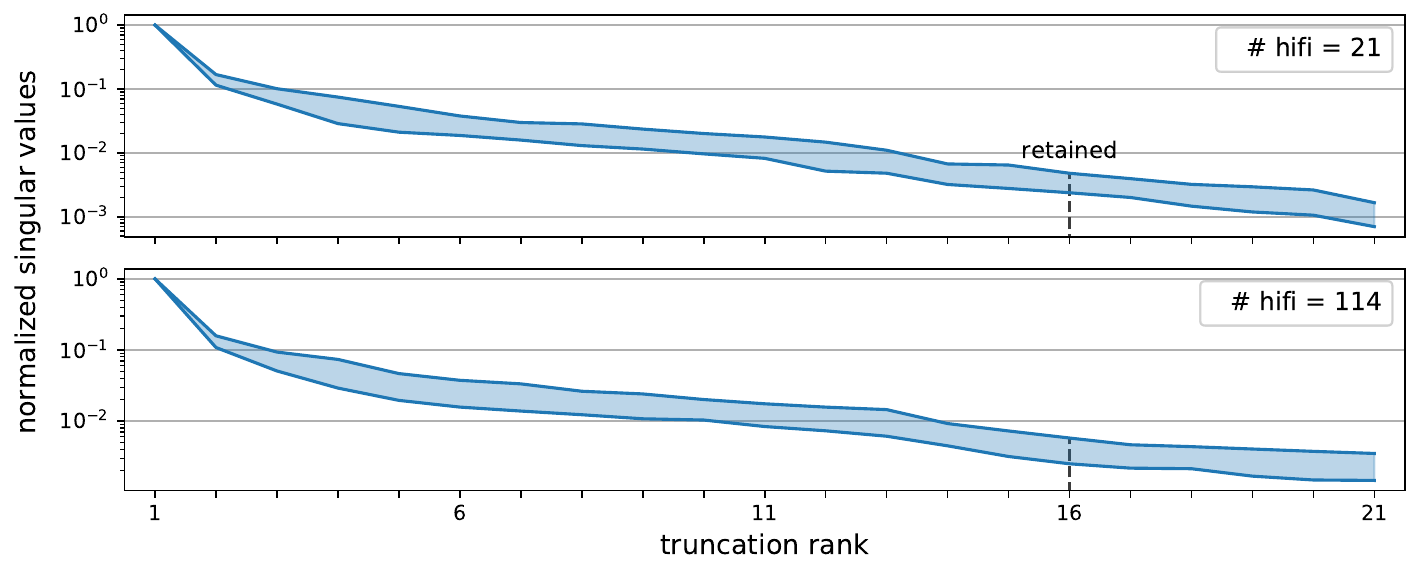}
	\caption{Evolution of the minimum and maximum for the normalized singular values decay of the stress tensor components, for the full ship with 10 parameters.
    Values are reported for the high-fidelity database before and after optimization.
    }
	\label{fig:fullship_10p_svalues_evo}
\end{figure}

\begin{figure}[hbt!]
	\centering
	\includegraphics[trim=0 0 0 0, clip, width=.97\textwidth]{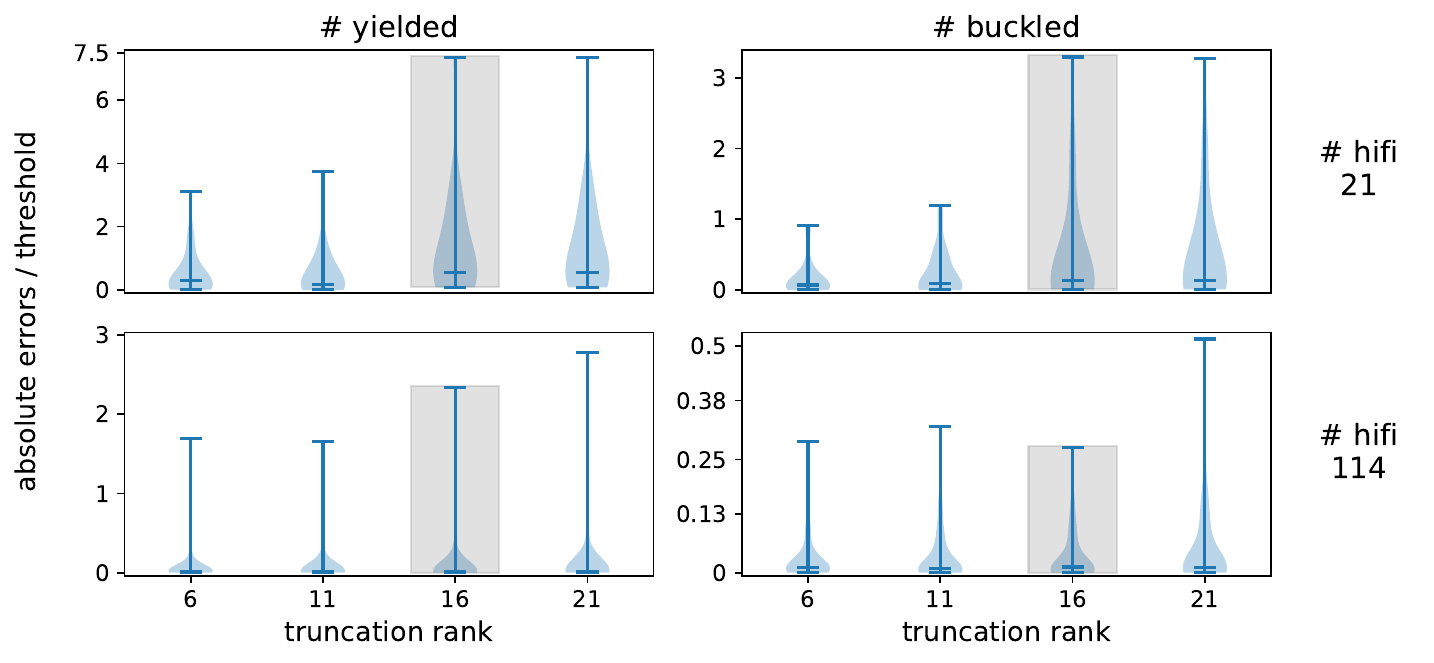}
	\caption{Evolution of the 5-fold cross-validation error distribution for the full ship with 10 parameters.
    Values are reported for the high-fidelity database before and after optimization. The gray box identifies the chosen truncation rank.
    }
	\label{fig:fullship_10p_errors_evo}
\end{figure}

Figure~\ref{fig:fullship_10p_svalues_evo} shows the evolution of the singular values decay for the model with the initial parameterization.
With the initial sampling of 21 configurations, the normalized singular values from rank 14 and above fall below the $10^{-2}$ threshold, but the truncation rank was set to 16 as a precaution against the possible complexity arising from the larger number of parameters.
As the optimizations explores the design space, this proved to be too conservative, and the singular values remained below the $10^{-2}$ threshold for ranks 14 and above.
Figure~\ref{fig:fullship_10p_errors_evo} shows the evolution of the errors on the yielded and buckled elements.
The maximum errors on the yielded elements are much larger than those measured on the buckled elements, but the addition of high fidelity configurations to the database determines a substantial concentration of the errors close to zero.
Unexpectedly, the error distributions corresponding to truncation rank 6 appear comparable to those obtained using higher truncation ranks.
This can be linked to the effect of the selection bias from the optimization procedures: 4 of the original 10 parameters showed very little influence on the QoIs, and neither multi- nor single-objective optimization moved them beyond their minimum values.
This results in the high-fidelity database underutilizing 4 out of 10 dimensions of the parametric domain, and as the number of validated configurations grows, the prediction errors of the ROMs reflect that the actual rank is about 4 positions lower than the number of parameters.

\begin{figure}[hbt!]
	\centering
	\includegraphics[trim=0 0 0 0, clip, width=.93\textwidth]{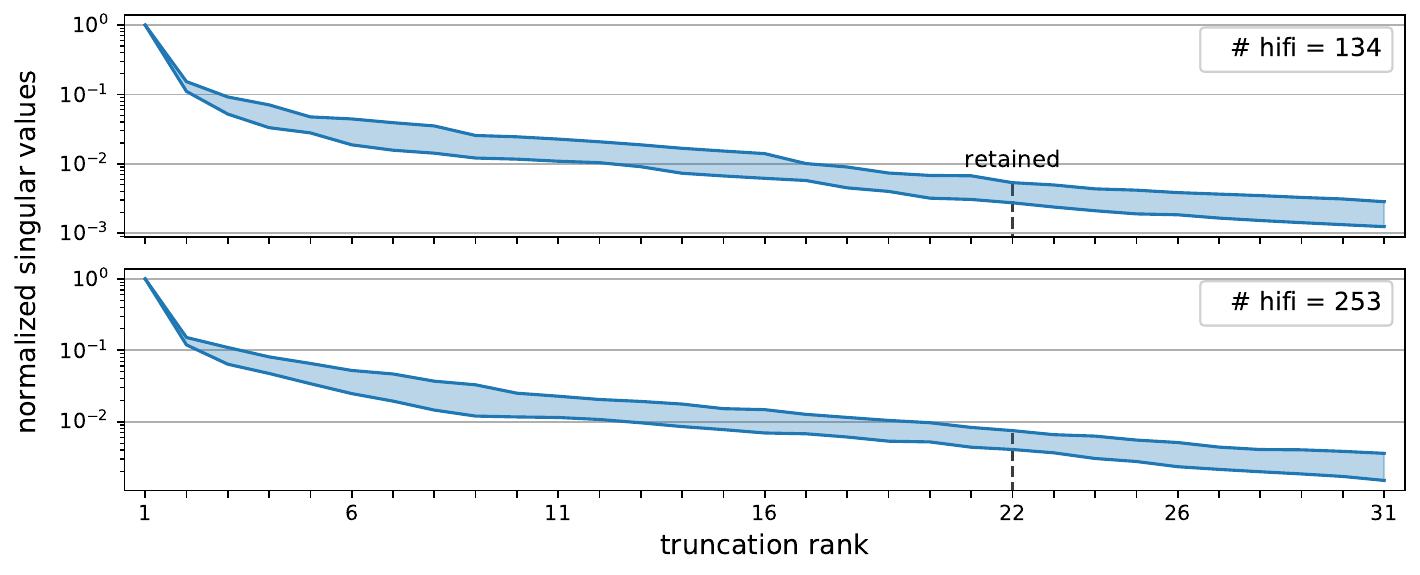}
	\caption{Evolution of the minimum and maximum for the normalized singular values decay of the stress tensor components, for the full ship with 16 parameters.
    Values are reported for the high-fidelity database before and after optimization.
    }
	\label{fig:fullship_10p_to_16p_svalues_evo}
\end{figure}

After the first reparameterization, the truncation rank is again increased by the number of new decision variables.
Of the 10 parameters, 6 were split in two, with the remaining 4 deemed not valuable for further refinement.
Again, the automatic adjustment of the truncation rank to 22 proves sufficient with respect to the $10^{-2}$ threshold on the normalized singular values.
This situation is depicted in Figure~\ref{fig:fullship_10p_to_16p_svalues_evo}.
As for the errors on the QoIs, the same trend as for the previous parameterization is seen, with the error distributions being very similar for different truncation ranks, but showing a concentration close to zero as the high-fidelity database size increases.

\begin{figure}[hbt!]
	\centering
	\includegraphics[trim=0 0 0 0, clip, width=.97\textwidth]{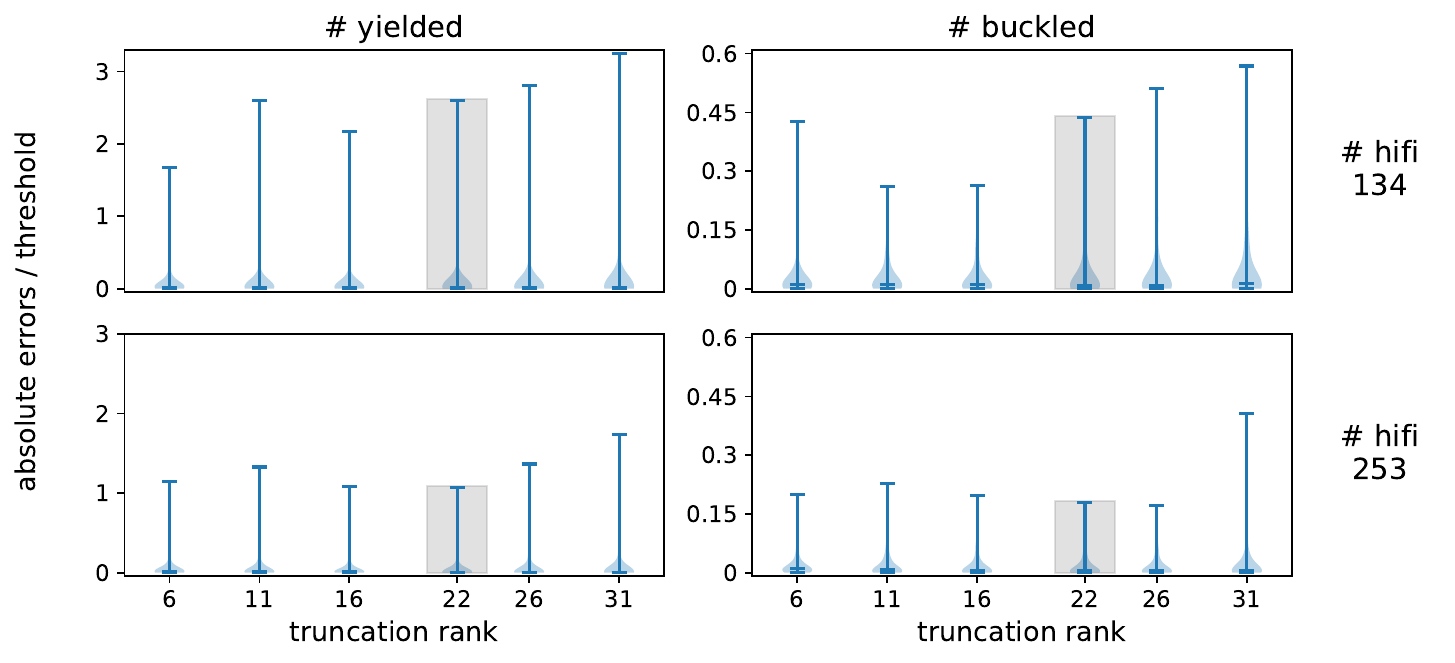}
	\caption{Evolution of the 5-fold cross-validation error distribution for the full ship with 16 parameters.
    Values are reported for the high-fidelity database before and after optimization. The gray box identifies the chosen truncation rank.
    }
	\label{fig:fullship_10p_to_16p_errors_evo}
\end{figure}

\end{document}